\documentclass[a4paper]{article}
\usepackage[latin1]{inputenc} %
\usepackage[T1]{fontenc} %
\usepackage{RR}
\usepackage{hyperref}

\usepackage{amsfonts,latexsym}
\usepackage{theorem}
\usepackage{amsmath}
\RRdate{Novembre 2009}

\RRauthor{
Laurent Baratchart\thanks[sfn]{INRIA, BP 93, 
06902 Sophia-Antipolis Cedex, FRANCE} 
  \and
Juliette Leblond\thanksref{sfn}
\and 
Fabien Seyfert\thanksref{sfn}
}
\authorhead{L. Baratchart \& J. Leblond \& F. Seyfert}
\RRtitle{Probl\`emes extr\'emaux contraints dans l'espace de Hardy $H^2$ et formules de Carleman}
\RRetitle{Constrained extremal problems in the Hardy space $H^2$ and Carleman's formulas}
\titlehead{Constrained extremal problems and Carleman's formulas} 
\RRresume{Nous \'etudions des probl\`emes d'approximation sur un strict sous-ensemble du cercle par des fonctions analytiques de l'espace de Hardy $H^2$ du disque unit\'e (de $\CC$), soumises en module à une contrainte ponctuelle sur la partie compl\'ementaire du cercle.
Des r\'esultats d'existence, d'unicit\'e et de saturation de la contrainte sont \'etablis, ainsi qu'une \'equation aux points critiques qui permet de proposer une formule de dualit\'e. Nous consid\'erons enfin une version polyn\^omiale (de dimension finie) de ces probl\`emes extr\'emaux born\'es.}
\RRabstract{We study  some approximation problems on a strict subset of the circle by analytic functions of the Hardy space $H^2$ of the unit disk (in $\CC$), whose modulus satisfy a pointwise constraint on the complentary part of the circle. 
Existence and 
uniqueness results, as well as pointwise saturation of the constraint, are established.
We also derive a critical point equation which gives rise to a dual formulation of 
the problem. We further
compute directional derivatives for this functional as a computational means to approach the
issue.
We then consider a
finite-dimensional polynomial version of the bounded extremal problem.}
\RRmotcle{Espaces de Hardy, fonctions analytiques, approximation, probl\`emes extr\'emaux born\'es, formules de Carleman.}
\RRkeyword{Hardy spaces, analytic functions, approximation, bounded extremal problems, Carleman formulas.}
\RRprojet{Apics}
\RRdomaine{3} 
\RRtheme{Modélisation, optimisation et contrôle de systèmes dynamiques}
\RCSophia 
\def\CC{{\mathbb C}}
\def\DD{{\mathbb D}}
\def\NN{{\mathbb N}}

\def\RR{{\mathbb R}}
\def\TT{{\mathbb T}}

\newtheorem{theorem}{Theorem}
\newtheorem{lemma}{Lemma}
\newtheorem{corollary}{Corollary}
\newtheorem{proposition}{Proposition}
\newcommand{\boite}{\mbox{} \hfill \mbox{\rule{2mm}{2mm}}}
\begin{document}
\RRNo{7087}
\makeRR   

\section{Introduction}
\label{sec1}
If $D$ is a finitely connected plane domain with rectifiable boundary
$\partial D$, a holomorphic function $f$ in the Smirnov class 
$\mathcal{E}^1(D)$ 
can be recovered from its boundary values by the Cauchy formula \cite{Duren}. 
When the boundary values are only known on a strict subset $I$ of 
$\partial D$ having positive linear measure, 
they still define $f$ uniquely but the recovery cannot be achieved 
in closed form. 
In fact, it becomes a special case of a classical ill-posed issue namely 
the Cauchy problem 
for the Laplace equation. This issue is quite important in
physics and engineering  \cite{Isakov,Lavrentiev,Par97}.

Following an original idea of Carleman, one approach to the recovery of $f$
from its knowledge on $I$ is to introduce an auxiliary ``quenching'' function
$\varphi$, holomorphic and bounded in $D$, 
such that $|\varphi|\equiv 1$ 
a.e. on $\partial D\setminus I$ and $|\varphi|>1$ in $D$; such a 
function is easily
constructed by solving a Dirichlet problem for $\log|\varphi|$.
In \cite{GoKry}, it was proven
by Goluzin and Krylov  that
\begin{equation}
\label{presCarleman}
f(z)=\lim_{n\to\infty}\,f_n(z),~~~~\mbox{where}~~\,f_n(z)
\stackrel{\Delta}{=}\frac{1}{2\pi}\int_I
\left(\frac{\varphi(\xi)}{\varphi(z)}\right)^n\frac{f(\xi)}{\xi-z}\,d\xi,
~~~~~~~~z\in D,
\end{equation}
the convergence being locally uniform in $D$. Cauchy integrals like those
defining $f_n$ in (\ref{presCarleman}) are called
\emph{Carleman's formulas} \cite{Aizenberg}.
On the unit disk $\DD$ where $\mathcal{E}^p(\DD)$ coincides with
the Hardy class $H^p$, see \cite{Duren},
it is proved in \cite{Patil} that if $f\in H^p$ with $1<p<\infty$,
then the convergence actually holds in $H^p$ .

Two questions arise naturally, namely what is the meaning of $f_n$
for \emph{fixed} $n$, and what is its asymptotic behaviour
if $f\in L^p(I)$ is \emph{not} the trace of a Hardy function?
On $\DD$, when $f\in L^2(I)$ and $\varphi$ is a quenching function
with the additional property that $|\varphi|$ is constant a.e. on $I$, 
it was proven in \cite{partII} that $f_n$ is closest to $f$ in $L^2(I)$-norm
among all $g\in H^2$ such that $\|g\|_{L^2(\TT\setminus I)}\leq 
\|f_n\|_{L^2(\TT\setminus I)}$, where $\TT$ denotes the unit circle.
In the present paper, among other things,
 we will see that if $\varphi$ is holomorphic
and bounded on $\DD$ together with its inverse, then $f_n$ is closest to
$f$ w.r.t. the weighted $L^2(|\varphi_{|_I}|^2,I)$-norm among all $g\in H^2$ such that 
$|g|\leq |f_n|$ a.e. on $\TT\setminus I$. These extremal properties of $f_n$
are all the more remarkable than  Carleman's formulas were originally 
defined without reference to optimization.

Still one point is unsatisfactory, namely the extremal properties
of $f_n$ we just mentioned are implicit in that
the level of the pointwise constraint on $\TT\setminus I$ is $|f_n|$ itself.
This is why we make a slight twist and we rather
investigate on $\DD$ the following extremal problem.
Let $I\subset\TT$ be a subset of positive Lebesgue measure 
and set  $J = \TT \setminus I$ for the complementary subset. 
The question we raise is the following.

\begin{description}
\item[$BEP$:]
Given $f \in L^2(I)$ and $M\in L^2(J)$, $M\geq 0$, 
find $g_0 \in H^2$ such that
 $|g_0(e^{i \theta})| \leq M(e^{i \theta})$
a.e. on $J$ and
\begin{equation} 
\label{problem}
\|f - g_0\|_{L^2(I)} = \min_{\stackrel{g \in H^2}{|g| \leq M \,
    \mbox{a.e. on} \, J}} \|f - g\|_{L^2(I)} \, .
\end{equation}  
\end{description}
This should be compared with the so-called \emph{bounded extremal 
problems} $BEP_p$.
studied in \cite{ablinria,partII,partI} for $1\leq p\leq\infty$.

\begin{description}
\item[$BEP_p$:]
Given $f \in L^p(I)$, $\psi\in L^p(J)$ and a positive 
constant $C$, find $g_0 \in H^p$ such that
\[\|g_0-\psi\|_{L^p(J)} \leq C~~ \mathrm{and}\]
\begin{equation} 
\label{problempq}
\|f - g_0\|_{L^p(I)} = \min_{\stackrel{g \in H^p}{\|g-\psi\|_{L^p(J)}
\leq C }} \|f - g\|_{L^p(I)} \, .
\end{equation}  
\end{description}

Note that in problem $BEP$ (\ref{problem}), we did not introduce a reference function
$\psi$ on $J$ as was done in $BEP_p$ (\ref{problempq}). 
While it is straightforward to handle such a generalization when
$\psi$ is the trace on $J$ of a $H^2$ function, the general case conceals further 
difficulties that are left here for further research.


When $I$ is of full measure, both problems (\ref{problem}) and
(\ref{problempq}) 
reduce to classical extremal problems, 
see {\it e.g.} \cite{Duren,Garnett}. Therefore we limit our discussion
to the case where $J$ has positive measure. 

The first reference dealing with
bounded extremal problems seems to be \cite{KN},
where $BEP_2$ is studied
for $f=0$ and $I$ an interval on the half-plane rather than the disk.
The case $\psi=0$ is solved in \cite{ablinria} using Toeplitz operators, and
error rates when $C$ goes large and $I$ is an arc
can be found in \cite{BLPprep}. Weighted versions of $BEP_2$ for $L^2(|\varphi|^2)$-norms as those discussed above were also solved in \cite{LO98}.

The general version $BEP_p$ in the range $1\leq p<\infty$ is taken up in \cite{partII}
where the link with Carleman's formulas is pointed out, while
existence and uniqueness results are also presented.
Reformulations of $BEP_p$
in abstract Hilbert or smooth Banach space settings were carried out
in \cite{CP2002,CPM2002,LP99,MSmith}, leading 
to the construction of backward minimal vectors and hyperinvariant subspaces 
for certain classes of operators that need not be compact nor quasinilpotent, 
thereby generalizing \cite{AE98}.
Versions of $BEP_2$ where the constraint bears on
the imaginary part rather than the modulus, useful among other things
to approach inverse Dirichlet-Neumann problems, are presented in \cite{imh2, LMP08}.
Together with  meromorphic
generalizations, problem $BEP_p$ was studied
in \cite{BS02} for $p \geq 2$, while problem $BEP_\infty$ was studied
in \cite{partI,aakblp}, with related completion issues. 

A major incentive to study 
$BEP_p$ came from engineering problems, more precisely
from questions pertaining to system identification and design.
This motivation is quite explicit in \cite{KN}, 
and all-pervasive in 
\cite{ablinria, BLPprep, partI, aakblp, thesefab} whose
results have been used effectively  to identify hyperfrequency filters 
\cite{RRbglosw}. The connection with identification is more transparent on 
the half plane, where $f$ represents the so-called transfer-function of a
linear dynamical system as measured in the frequency bandwidth $I$
using harmonic identification techniques. Recall 
that a linear dynamical system is just a convolution operator, and 
that its transfer function is the Fourier-Laplace transform of its kernel
\cite{DFT}. Now, by the Paley-Wiener and 
Hausdorff-Young theorems, causality and $L^r\to L^s$ stability
of the system cause $f$ to belong to the Hardy class $\mathcal{H}^p$
of the half plane with $1/p=1/r-1/s$, as soon as the latter is
less than or equal to $1/2$. Because $f$ is only known 
up to modelling and measurement errors, one is led to
approximate the data on $I$ by a $\mathcal{H}^p$ function 
while controlling its deviation from some reference behaviour $\psi$ 
outside $I$, which is precisely the analog of
(\ref{problempq}) on the half-plane.
It is mapped to $BEP_p$ {\it via} the isometry 
$g \mapsto (1+w)^{-2/p}
g((w-1)/(w+1))$ from $H^p$ onto ${\cal H}^p$. More on the relations between
Hardy spaces, system identification and control can be found in
\cite{Fuhrmann, Nikolskii2,Par97}.
Note that in  $BEP_p$, it is indeed essential 
to bound the behaviour of $g_0$ on $J$, for
traces of Hardy functions are dense in $L^p(I)$ (in $C(I)$ if
$p=\infty$) so that
$BEP_p$ has no solution 
if $C=\infty$ unless $f$ is already the trace of a Hardy function.
In practice, since modelling and measurement errors will prevent this 
from ever happening, the error $\|f-g\|_{L^p(I)}$ can 
be made arbitrarily small at the cost of $\|g\|_{L^p(J)}$ becoming
arbitrarily large, which does not make for a valid identification scheme.

The present paper 
deals with a mixed situation, 
where an integral criterion is minimized on $I$ under a pointwise
constraint on $J$. 
Here again, the motivation of the authors
stems from system identification. Indeed, the $L^2$ norm 
on $I$ has a probabilistic  interpretation, being the variance
of the output when the system is fed by 
noise whose spectrum is uniformly distributed in the bandwidth, that suits
some classical framework. On another hand,
it is often the case that the transfer function has to meet uniform bounds 
for physical reasons. For instance when identifying
a passive device, it should be less than 1 at all frequencies.
This way one is led to consider problem  (\ref{problem}) with $M\equiv1$ (which is $BEP_{2, \infty}$ below).

Such issues and motivations are also the topics of the recent work \cite{schneck1, schneck2}, where  $BEP$-like  problems are considered in $L^p(I)$, for $1 \leq p \leq \infty$, with a pointwise constraint acting on the whole $\TT$, when the approximated function $f$ and the constraint $M$ are assumed to be continuous functions (on $\TT$ and $I$, respectively), with $M >0$ and $|f| \leq M$ on $I$.

Problem (\ref{problem})
is considerably more difficult to analyze than $BEP_2$,
due to the fact that pointwise evaluation is not smooth 
--actually not even defined-- in $L^2(J)$. Its solution depends in a 
rather deep fashion on the multiplicative structure of Hardy functions and
all our results, beyond existence and uniqueness, will hold
under the extra-assumption that the boundary of 
$I$ has measure zero. We do not know the extend to which this assumption can 
be relaxed.\\

The organization of the paper is as follows. In section \ref{notprel} we set
up some notation and recall standard properties of $H^p$-spaces, $BMO$
and conjugate functions. Section \ref{wposed} deals with existence and 
uniqueness issues, as well as pointwise saturation of the constraint.
In section \ref{critpsec} we establish an analog, 
in this nonsmooth and infinite-dimensional context, of the familiar
critical point equation from convex analysis. It gives rise to a saddle-point 
characterization of the optimal value that yields a dual formulation of 
the problem. 
The latter is connected in section \ref{Carlsec} to Toepliz operators
and Carleman's formulas, and
used to compute the concave dual functional
whose maximization is tantamount to solve the problem. We further
compute directional derivatives for this functional as a means to approach the
issue from a computational point of view. Section \ref{const} is devoted to a
finite-dimensional polynomial version of (\ref{problem}), valid when $I$ is a 
union of arcs, which is of interest 
in its own right and provides an alternative way to constructively approximate
the solution to the original problem.
\newpage
\section{Notations and preliminaries}
\label{notprel}
Let $\TT$ be the unit circle endowed with the
normalized Lebesgue measure $\ell$, 
and $I$ a subset of $\TT$ such that $\ell(I)>0$ 
with complementary subset $J=\TT\setminus I$. 
To avoid dealing with trivial instances of problem (\ref{problem}) 
\emph{ we assume throughout that $\ell(J)>0$.}

If $h_1$ (resp. $h_2$) is a function defined on a set containing $I$ 
(resp. $J$), we use the notation $h_1\vee h_2$ 
for the concatenated function, defined on the whole of $\TT$, which is $h_1$
on $I$ and $h_2$ on $J$. 

For $E\subset\TT$, we let $\partial E$ and 
$\stackrel{\circ}{E}$ denote respectively the boundary and the interior
of $E$ when viewed as a subset of $\TT$; we also write $\chi_E$ for the
characteristic function of $E$ and
$h_{|_E}$ to mean the restriction to $E$ of a function $h$ defined on a set
containing $E$. 

When $1\leq p\leq\infty$, we write $L^p(E)$ for the familiar Lebesgue space of
(equivalence classes of a.e. coinciding) complex-valued measurable functions 
on $E$ with 
finite $L^p$ norm, and we indicate by $L^p_{\RR}(E)$ the real subspace of
real-valued functions.
Likewise $C(E)$ stands for the space of 
complex-valued continuous functions on
$E$, while $C_\RR(E)$ indicates real-valued 
continuous functions. The norm on $L^p(E)$ is denoted by
$\|~\|_{L^p(E)}$, and if $h$ is defined on a set
containing $E$ we write for simplicity $\|h\|_{L^p(E)}$ to
mean $\|h_{|_E}\|_{L^p(E)}$. When $E$ is compact
the norm of $C(E)$ is the {\it sup} norm.

Recall that the Hardy space $H^p$ is the closed
subspace of $L^p(\TT)$ consisting of functions whose Fourier coefficients
of strictly negative index do vanish. These are the nontangential limits of 
functions analytic in the unit disk $\DD$ having uniformly bounded $L^p$ 
means over all circles centered at $0$ of radius less than 1. 
The correspondence is one-to-one and, using this 
identification, we alternatively regard members of $H^p$ as holomorphic
functions in the variable $z\in\DD$. 
This extension is obtained from the values on $\TT$
through a Cauchy as well as a Poisson integral \cite[ch. 17, thm 11]{Rudin}, 
namely if $g\in H^p$ then, for $z\in\DD$:
\begin{equation}
\label{Cauchy1}
g(z) = \frac{1}{2\,i \, \pi} \, \int_{\TT}
\frac{g(\xi)}{\xi-z}\,d\xi \mbox{ and }
g(z)=\frac{1}{2 \, \pi} \, \int_{\TT} \mbox{\rm Re}\left\{\frac{e^{i\theta} + 
z}{e^{i\theta} - z} \right\}\, g(e^{i\theta})\, d\theta \, .
\end{equation}
Because of this Poisson representation, $g(re^{i\theta})$ converges to $g(e^{i\theta})$ in 
$L^p(\TT)$ as soon as $1\leq p<\infty$. Moreover, (\ref{Cauchy1}) entails
that, for $1\leq p\leq\infty$, a Hardy function $g$ is uniquely determined, 
up to a purely imaginary constant, by its real part $h$ on $\TT$:
\begin{equation}
g(z) = i\mbox{\rm Im}g(0)\,+\,\frac{1}{2 \, \pi} \, \int_{\TT} \frac{e^{i\theta} + 
z}{e^{i\theta} - z} \, h(e^{i\theta})\, d\theta \, ,~~ ~~z\in\DD.
\label{RHt}
\end{equation}
The integral in the right-hand side of (\ref{RHt}) is called the
\emph{Riesz-Herglotz transform} of $h$ and, whenever $h\in L_{\RR}^1(\TT)$, it
defines a holomorphic function in $\DD$ which is real at 0 and whose
nontangential limit exists a.e. on $\TT$ with real part equal to $h$.
However, only if  $1<p<\infty$ is it guaranteed  that $g\in H^p$
when $h\in L_{\RR}^p(\TT)$. In fact,
the Riesz-Herglotz transform assumes the form 
$h(e^{i\theta})+i\widetilde{h}(e^{i\theta})$ a.e. on $\TT$, where the
real-valued function $\widetilde{h}$ is said to be \emph{conjugate} to $h$,
and the property that $\widetilde{h}\in L_{\RR}^p(\TT)$ whenever
$h\in L_{\RR}^p(\TT)$ holds true for $1<p<\infty$ but not 
for $p=1$ nor $p=\infty$. The map $h\to\widetilde{h}$ is called the 
\emph{conjugation operator}, and for $1<p<\infty$  it is bounded 
$L_{\RR}^p(\TT)\to L_{\RR}^p(\TT)$ by a theorem of M. Riesz
\cite[chap. III, thm 2.3]{Garnett}; 
in this range of exponents, we will denote its norm
by $K_p$. It follows easily from Parseval's relation that $K_2=1$, but it is rather subtle that
$K_p=\tan(\pi/(2p))$ for $1<p\leq2$ while $K_p=\cot(\pi/(2p))$ for $2\leq p<\infty$
\cite{Pichorides}.

A sufficient condition for $\widetilde{h}$ to be in $L^1(\TT)$ is that $h$ belongs to
the the so-called \emph{Zygmund class} $L\log^+ L$, consisting of measurable functions
$\phi$ such that $\phi\log^+|\phi|\in L^1(\TT)$ where we put $\log^+ t=\log t$ if $t\geq1$ 
and 0 otherwise. More precisely, if we denote by $m_h$ the \emph{distribution function} of $h$
defined on $\RR^+$ with values in $[0,1]$ according to the formula
\[m_h(\tau)=\ell\left(\{\xi\in\TT;~|h(\xi)|>\tau\}\right),\]
and if we further introduce the non-increasing rearrangement of $h$ given by
\[h^*(t)=\inf\{\tau;~m_h(\tau)\leq t\},~~~~~~t\geq0,\]
it turns out that $h\in L\log^+ L$ if and only if the quantity
\begin{equation}
\label{defLlogL}
\|h\|_{L\log^+L}\stackrel{\Delta}{=}\int_0^1h^*(t)\log(1/t)\,dt
\end{equation} 
is finite \cite[lem. 6.2.]{BeSh}, which makes $ L\log^+ L$ into a Banach function space.
Then, it is a theorem of Zygmund \cite[cor. 6.9.]{BeSh} that
\begin{equation}
\label{LlogLconj}
\|\widetilde{h}\|_{L^1(\TT)}\leq C_0\|h\|_{L\log^+L}
\end{equation}
for some universal constant $C_0$. A partial converse, due to M. Riesz, asserts that
if a real-valued $h$ is bounded from below and if moreover $\widetilde{h}\in L^1(\TT)$,
then $h\in L\log^+ L$ \cite[cor. 6.10]{BeSh}.

We mentioned already that
$\widetilde{h}$ needs not be bounded if $h\in L^\infty_{\RR}(\TT)$. In this
case all one can say in general is that $\widetilde{h}$ has \emph{bounded mean
oscillation}, meaning that $\widetilde{h}\in L^1(\TT)$ and
\[\|\widetilde{h}\|_{BMO}\stackrel{\Delta}{=}
\sup_{E}\frac{1}{\ell(E)}\int_E|\widetilde{h}-\widetilde{h}_E|\,d\theta<\infty,
~~~~~~{\rm with}~~\widetilde{h}_E\stackrel{\Delta}{=}\frac{1}{\ell(E)}
\int_E\widetilde{h}\,d\theta,\]
where the {\it supremum} is taken over all subarcs $E\subset\TT$.
Actually \cite[chap. VI, thm 1.5]{Garnett}, there is a universal 
constant $C_1$ such that
\[
\|\widetilde{h}\|_{BMO}\leq C_1\|h\|_{L^\infty(\TT)}.
\]
The subspace of $L^1(\TT)$ consisting of functions whose $BMO$-norm is finite is called
$BMO$ for short. Notice that $\|~\|_{BMO}$ is a genuine norm 
modulo additive constants
only. A theorem of F. John and L. Nirenberg \cite[ch. VI, thm. 2.1]{Garnett}
asserts there are positive constants $C$, $c$, such that, 
for each real-valued $\varphi\in BMO$, every arc $E\subset\TT$, and any $x>0$,
\begin{equation}
\label{JN}
\frac{\ell\left(\{t\in E:~|\varphi(t)-\varphi_E|>x\}  \right)}{\ell(E)}\leq 
C\exp\left(\frac{-cx}{\|\varphi\|_{BMO}}\right).
\end{equation}
Conversely, if (\ref{JN}) holds for some finite $A>0$ in place of $\|\varphi\|_{BMO}$, 
every arc $E$ and any $x>0$, then $\varphi\in BMO$ and $A\sim \|\varphi\|_{BMO}$. 
The John-Nirenberg theorem easily implies that $BMO\subset L^p$ for
all $p<\infty$.
The space of $H^1$-functions whose boundary values lie 
in $BMO$ will be denoted by $BMOA$, and $BMOA/\CC$ is a Banach space equipped
with the $BMO$-norm. Clearly $BMOA\subset H^p$ for $1\leq p<\infty$, and
$h+i\widetilde{h}\in BMOA$ whenever $h\in L^\infty(\TT)$. 
A sufficient condition for the
boundedness of  $\widetilde{h}$ is that $h$ be Dini-continuous;
recall that a function $h$ defined on
$\TT$ is said to be Dini-continuous if $\omega_h(t)/t \in L^1([0,\pi])$, where
\[\omega_h(t)=\sup_{|\theta_1-\theta_2|\leq t}~ 
\left|h\bigl(e^{i\theta_1}\bigr)-
h\bigl(e^{i\theta_2}\bigr)\right|,~~~~t\in[0,\pi],\]
is the modulus of continuity of $h$. Specifically \cite[chap. III, thm
1.3]{Garnett}, it holds that 
\begin{equation}
\label{modDini}
\omega_{\widetilde{h}}(\rho)\leq C_2\left(\int_0^\rho
  \frac{\omega_h(t)}{t}\,dt\,+\,
  \rho\int_\rho^\pi\frac{\omega_h(t)}{t^2}\,dt\right)
\end{equation}
where $C_2$ is a constant independent of $f$. From (\ref{modDini})
it follows easily that $\widetilde{h}$ is continuous if $h$ is
Dini-continuous, and moreover that
\begin{equation}
\label{majDini}
\|\widetilde{h}\|_{L^\infty(\TT)}\leq 
\omega_{\widetilde{h}}(\pi)\leq C_2\int_0^\pi \frac{\omega_h(t)}{t}\,dt,
\end{equation}
where the first inequality comes from the fact that $\widetilde{h}$ is
continuous on $\TT$ and therefore vanishes at some point 
since it has zero-mean.

We now turn to multiplicative properties of Hardy functions.
It is well-known 
(see {\it e.g.} \cite{Duren,Garnett,Koosis}) 
that a nonzero $f\in H^p$ can be uniquely
factored as $f=jw$ where
\begin{equation}
\label{def-outer}
w(z)=\exp\left\{\frac{1}{2\pi}\int_0^{2\pi}\frac{e ^{i\theta}+z}{e
    ^{i\theta}-z}\log|f(e ^{i\theta})|\, d\theta\right\}
\end{equation}
belongs to $H^p$ and is called the {\em outer factor} of $f$, while
$j\in H^\infty$ has modulus 1 a.e. on $\TT$ and is called the
{\em inner factor} of $f$. The latter may be further decomposed as
$j=bS_{\mu}$, where
\begin{equation}
\label{def-Blaschke}
b(z)=e^{i\theta_0}
z^k\prod_{z_l\neq0}\frac{{-\bar z}_l}{|z_l|}\,\frac{z-z_l}{1-{\bar z}_lz}
\end{equation}
is the \emph{Blaschke product}, with order $k\geq0$ at the origin,
associated to the sequence $z_l\in\DD\setminus\{0\}$
and to the constant $e^{i\theta_0}\in{\bf T}$, while
\begin{equation}
\label{def-singulier}
S_{\mu}(z)=\exp\left\{-\frac{1}{2\pi}\int_0^{2\pi}\frac{e ^{i\theta}+z}{e
    ^{i\theta}-z}\, d\mu(\theta)\right\}
\end{equation}
is the {\em singular inner factor} associated with $\mu$, a positive measure
on $\TT$ which is singular with respect to Lebesgue measure. 
The $z_l$ are of course the zeros of $f$ in $\DD\setminus\{0\}$, counted
with their multiplicities, while $k$ is the order of the zero at 0. 
If there are infinitely many zeros,
the convergence of the product $b(z)$ in $\DD$
is ensured by the condition
$\sum_l(1-|z_l|)<\infty$ which holds automatically when 
$f\in H^p\setminus\{0\}$. If there are only finitely many $z_l$, 
we say that (\ref{def-Blaschke}) is a finite Blaschke product; note that
a finite Blaschke product may alternatively be defined as a rational function
of the form $q/q^R$, where $q$ is an algebraic polynomial 
whose roots lie in $\DD$ and $q^R$ 
indicates the \emph{reciprocal polynomial} given by 
$q^R(z)=z^n\overline{q(1/{\bar z})}$ if $n$ is the degree of $q$. 
The integer $n$ is also called the degree of the Blaschke product.

That $w(z)$ in (\ref{def-outer}) is well-defined rests on the fact that
$\log|f|\in L^1$ if $f\in H^1\setminus\{0\}$; this also
entails that a $H^p$ function cannot vanish 
on a subset of strictly positive Lebesgue measure on $\TT$ unless it is 
identically 
zero. For simplicity, we often say that a function is outer (resp. inner) 
if it is equal to its outer (resp. inner) factor.

Intimately related to Hardy functions is the Nevanlinna class $N^+$ consisting
of holomorphic functions in $\DD$
that can be factored as $jE$, where $j$ is an inner function and $E$ an outer 
function of the form
\begin{equation}
\label{outa}
E(z)=\exp\left\{\frac{1}{2\pi}\int_0^{2\pi}\frac{e ^{i\theta}+z}{e
    ^{i\theta}-z}\log \rho(e ^{i\theta})\, d\theta\right\}\,,
\end{equation}
with $\rho$ a positive function
such that $\log\rho\in L^1(\TT)$ (although $\rho$ itself need not
be summable). Such a function again has nontangential limits of modulud $\rho$ a.e. on $\TT$
that serve as definition of its boundary values. The Nevanlinna class
will be instrumental to us in that
$N^+\cap L^p(\TT)=H^p$, see  for example \cite[thm 2.11]{Duren} or 
\cite[5.8, ch.II]{Garnett}. Thus formula (\ref{outa}) defines a $H^p$-function
if, and only if $\rho\in L^p(\TT)$.
A useful consequence is that, whenever 
$g_1\in H^{p_1}$ and $g_2\in H^{p_2}$, we have $g_1g_2\in H^{p_3}$ if,
and only if $g_1g_2\in L^{p_3}$. In particular 
$g_1g_2\in H^{p_3}$ if $1/p_1+1/p_2=1/p_3$.

It is a classical fact \cite[ch. II, sec. 1]{Garnett}
that a function $f$ holomorphic in the unit disk
belongs to $H^p$ if, and only $|f|^p$, which is subharmonic in $\DD$,
has a harmonic majorant there.
This makes for a conformally invariant definition of Hardy spaces 
over general domains in $\overline{\CC}$. In this connection,
the Hardy space ${\bar H}^p$ of $\overline{\CC}\setminus\DD$ 
can be given a treatment parallel to $H^p$ using the
conformal map $z\mapsto1/z$. Specifically,
${\bar H}^p$  consists of $L^p$ functions
whose Fourier coefficients of strictly positive index do vanish; these
are, a.e. on $\TT$, the complex conjugates of $H^p$-functions, and
they can also be viewed as nontangential limits of 
functions analytic in $\overline{\CC}\setminus \overline{\DD}$
having uniformly bounded $L^p$ means 
over all circles centered at $0$ of radius bigger than $1$.
We also set $\overline{BMOA}={\bar H}^1\cap BMO$.
We further single out the subspace ${\bar H}_0^p$ of ${\bar H}^p$,
consisting of functions vanishing at infinity or, equivalently, 
having vanishing mean on $\TT$. Thus, a function belongs to  ${\bar H}_0^p$
if, and only if, it is a.e. on $\TT$ of the form 
$e^{-i\theta}\overline{g(e^{i\theta})}$ for some $g\in H^p$. For
$G\in{\bar H}^p_0$, the Cauchy formula assumes the form:
\begin{equation}
\label{Cauchy2}
G(z) = 
\frac{1}{2\,i \, \pi} \, \int_{\TT}
\frac{G(\xi)}{z-\xi}\,d\xi\,,~~~~z\in\overline{\CC}\setminus\overline{\DD}.
\end{equation}

If $E$ is a measurable subset of $\TT$, we set
\begin{equation}
\label{defdual}
<f,g>_E= \frac{1}{2\pi}\int_E f(e^{i\theta})\overline{g(e^{i\theta})}\,
d\theta
\end{equation}
whenever $f\in L^p(E)$ and $g\in L^q(E)$ with $1/p+1/q=1$. If $f$ and
$g$ are defined on a set containing $E$, we often write for simplicity
$<f,g>_E$ to mean $<f_{|_E},g_{|_E}>$.

The duality product $<~,~>_{\TT}$ makes $H^p$ and ${\bar H}_0^q$ 
orthogonal to each other, and reduces to the familiar scalar product 
on $L^2(\TT)\times L^2(\TT)$. We note in particular
the orthogonal decomposition:
\begin{equation}
\label{orthogdec}
L^2(\TT)=H^2\oplus{\bar H}_0^2.
\end{equation}
For $f\in C(\TT)$ and $\nu\in\mathcal{M}$, the space of complex 
Borel measures on 
$\TT$, we set 
\begin{equation}
\label{pairdual}
\nu{\bf .} f=\int_\TT f(e^{i\theta})\,d\nu(\theta)
\end{equation}
and this pairing induces an isometric 
isomorphism between $\mathcal{M}$ (endowed with the norm of the
total variation) and 
the dual of $C(\TT)$ \cite[thm 6.19]{Rudin}. If we let
$\mathcal{A}\subset H^\infty$ designate the disk algebra of
functions analytic in $\DD$ and continuous on $\overline{\DD}$, 
and if $\mathcal{A}_0$ indicates those functions in $\mathcal{A}$
vanishing at zero, it is easy to
see that $\mathcal{A}_0$ is the orthogonal space under
(\ref{pairdual}) to those measures whose Fourier
coefficients of strictly negative index do vanish. 
Now, it is a fundamental theorem 
of F. and M. Riesz that such measures have the form 
$d\nu(\theta)=g(e^{i\theta})\,d\theta$ 
with $g\in H^1$, so the Hahn-Banach theorem 
implies that $H^1$ is dual {\it via} (\ref{pairdual}) to the quotient space
$C(\TT)/\mathcal{A}_0$ \cite[chap. IV, sec. 1]{Garnett}. 
Equivalently, ${\bar H}^1_0$ is dual to $C(\TT)/\overline{\mathcal{A}}$ under the pairing
arising from the line integral:
\begin{equation}
\label{pairbar}
(\dot{f},F)=
 \frac{1}{2i\pi}\int_{\TT} f(\xi)F(\xi)\,d\xi\,,
\end{equation} 
where $F$ belongs to ${\bar H}^1_0$ and $\dot{f}$
indicates the equivalence class of $f\in C(\TT)$
modulo $\overline{\mathcal{A}}$. 
This entails that, contrary to $L^1(\TT)$, the spaces
$H^1$ and ${\bar H}^1_0$ enjoy a weak-* compactness property
of their unit ball.

Finally, we define the analytic and anti-analytic projections
${\bf P}_+$ and ${\bf P}_-$ on Fourier series by:
\[
{\bf P}_+\left(\sum_{n=-\infty}^{\infty}a_ne^{in\theta}\right)=
\sum_{n=0}^{\infty}a_ne^{in\theta},~~~~
{\bf P}_-\left(\sum_{n=-\infty}^{\infty}a_ne^{in\theta}\right)=
\sum_{n=-\infty}^{-1}a_ne^{in\theta}\,.
\]
Equivalent to the M. Riesz theorem is 
the fact that ${\bf P}_+:L^p\to H^p$ and ${\bf P}_-:L^p\to \bar{H}_0^p$
are bounded for $1<p<\infty$, in which
case they coincide with the Cauchy projections:
\begin{equation}
\label{Cauchyproj}
{\bf P}_+(h)(z)=\frac{1}{2i\pi}\int_{\TT}\frac{h(\xi)}{\xi-z}\,d\xi, 
{\bf P}_-(h)(s)=\frac{1}{2i\pi}\int_{\TT}\frac{h(\xi)}{s-\xi}\,d\xi, 
\end{equation}
for $z\in\DD$, $s \in\overline{\CC}\setminus\overline{\DD}$. 
When restricted to $L^2(\TT)$, the projections ${\bf P}_+$ and
${\bf P}_-$ are just the
orthogonal projections onto $H^2$ and $\bar{H}_0^2$ respectively. 
Likewise ${\bf P}_+:L^\infty\to BMOA$ and ${\bf P}_-:L^\infty\to \overline{BMOA}$
are also bounded.

Although ${\bf P}_\pm(h)$ needs not be the Fourier series of a function
when $h$ is merely in $L^1(\TT)$, it is nevertheless Abel summable 
almost everywhere to a function lying in $L^s(\TT)$
for $0<s<1$, and it can
still be interpreted as the trace of an analytic function in some
Hardy space of exponent $s$ that we did not introduce 
\cite[cor. to thm 3.2]{Duren}. To us it will be sufficient, 
when $h\in L^1$, to regard  ${\bf P}_\pm(f)$ as the Fourier 
series of a distribution. Finally, we record for further reference
the following elementary fact:
\begin{lemma}
\label{lemme_integral}
Let $v \in L^{1}(J)$ be such that ${\bf P}_+(0 \vee v)\in L^2(\TT)$.
Then, whenever $g\in H^2$ is such that $g\in L^2(I)\vee L^\infty(J)$,
it holds that 
\[<{\bf P}_+(0 \vee v)\,,\,g>_{\TT}\,=\,<v\,,\,g>_{J}. 
\]
\end{lemma}
{\sl Proof.}
Since by hypothesis ${\bf P}_+(0 \vee v)$ is a \emph{square summable function}
on $\TT$, we can define a function $u\in L^1(\TT)$ by the formula:
\[u=(0 \vee v) -{\bf P}_+(0 \vee v)\,,\]
and by the very definition of $u$ all its Fourier coefficients
of non-negative index do vanish hence  $u\in {\bar H}_0^1$.
In addition it is clear that $u_{|_I}\in L^2(I)$ and consequently,
if $g\in H^2$ is such that $g\in L^2(I)\vee L^\infty(J)$,
we have upon checking summability on $I$ and $J$ separately
that $u\overline{g}\in  {\bar H}_0^1$.
Therefore we get:
\[
<v\,,\,g>_{J}\,\,=\,\,<v\overline{g}\,,\,1>_{J}\,\,=\,\,
<(0\vee v)\,\overline{g}\,,\,1>_\TT
\,\,
\]
\[
=\,\,<u\,\overline{g}\,,\,1>_{\TT}+<{\bf P}_+(0 \vee v
)\,\overline{g}\,,\,1>_{\TT}
\]
\[
=\,\,<{\bf P}_+(0 \vee v )\,\overline{g}\,,\,1>_{\TT}
                      \,\,=\,\,<{\bf P}_+(0 \vee v)\,,\,g>_{\TT}
\]
where the next-to-last equality uses that the mean of the 
${\bar H}_0^1$-function $u\overline{g}$ is zero.
\hfill \boite\\  

\newpage
\section{Well-posedness of the bounded extremal problem $BEP$}
\label{wposed}
We first reduce problem $BEP$ (\ref{problem}) to a standard form $BEP_{2, \infty}$ where $M\equiv1$.
As the $\log$-modulus of a nonzero Hardy function is
integrable, we will safely assume that $\log M\in L^1(J)$
for otherwise the zero function is the only candidate approximant. Then, 
letting $w_M$ be the outer function with modulus $1$ on $I$
and $M$ on $J$, we have that $g$ belongs to $H^2$ and satisfies
$|g|\leq M$ a.e. on $J$ if, and only if $g/w_M$ lies in $H^2$
and satisfies $g/w_M\leq 1$ a.e. on $J$; it is so because
$g/w_M$ lies by construction in the Nevanlinna class $N^+$ whose 
intersection with $L^2(\TT)$ is $H^2$.
Altogether, upon replacing $f$ by $f/w_M$ and $g$ by $g/w_M$, we see 
that Problem (\ref{problem}) is equivalent to the 
following normalized  case which is the one we shall really work with.

\begin{description}
\item[$BEP_{2, \infty}$:] Given $f \in L^2(I)$, find $g_0 \in H^2$ such that 
$|g_0(e^{i \theta})| \leq 1$ a.e. on $J$ and
\begin{equation} 
\|f - g_0\|_{L^2(I)} = \min_{\stackrel{g \in H^2}{|g| \leq 1 \,
    \mbox{a.e. on} \, J}} \|f - g\|_{L^2(I)} \, .
\label{eq-approx}
\end{equation}  
\end{description}
Let us begin with a basic existence and uniqueness result:
\begin{theorem}
\label{existence}
Problem $BEP_{2, \infty}$ {\rm (\ref{eq-approx})} has a unique solution
$g_0$, and necessarily $\|g_0\|_{L^2(I)}\leq\|f\|_{L^2(I)}$. 
Moreover $\|g_0\|_{L^\infty(J)}=1$ unless
$f=g_{|_I}$ for some $g\in H^2$ such that 
$\|g\|_{L^\infty(J)}<1$.
\end{theorem}
\begin{corollary}
Problem $BEP$ {\rm (\ref{problem})} has a unique solution.
\end{corollary}
{\sl Proof of \rm{ Theorem \ref{existence}}.}
Define a convex subset of $L^2(I)$ by putting
${\cal C}:=\{g_{|_I};~~g\in H^2, ~\|g\|_{L^\infty(J)}\leq1\}$.
\emph{We claim} that ${\cal C}$ is closed. Indeed, 
let $\{g_n\}$ be a sequence in $H^2$, with $\|g_n\|_{L^\infty(J)}\leq1$,
that converges in $L^2(I)$ to some $\phi$.
Clearly $\{g_{n}\}$ is bounded in 
$L^2(\TT)$, therefore some subsequence $g_{k_n}$ converges weakly to
$g\in H^2$. Since $|{g_{k_n}}|\leq1$ on $J$,
we may assume upon refining the subsequence further
that it converges weak-* in $L^\infty(J)$ to a limit which can be none but
${g}_{|_J}$. By  weak-* compactness of balls in $L^\infty(J)$,
we get $\|g\|_{L^\infty(J)} \leq 1$,
hence ${g}_{|_I}\in{\cal C}$.
But ${g_{k_n}}_{|_I}$ {\it a fortiori}  converges weakly to ${g}_{|_I}$ in $L^2(I)$,
thus $\phi={g}_{|_I}\in{\cal C}$ \emph{as claimed}. 
By standard properties of the projection on a non-empty 
(for $0\in{\cal C}$) closed convex set in a Hilbert space,
we now deduce that the solution $g_0$ to (\ref{eq-approx}) uniquely exists, and 
is characterized by the variational inequality \cite[thm V.2.]{Brezis}
\begin{equation}
\label{varineg2}
{g_0}_{|_I}\in{\cal C}~~~~{\rm and}~~~~{\rm Re}\,
<f-g_0\,,\,\phi-g_0>_I\,\leq\,0,~~~~\forall \phi\in{\cal C}.
\end{equation}
Using $\phi=0$ in (\ref{varineg2}) and applying the Schwarz inequality yields $\|g_0\|_{L^2(I)}
\leq \|f\|_{L^2(I)}$.

Assume finally that $\|g_0\|_{L^\infty(J)}<1$. Given $h\in H^\infty$, 
$g_0+t h$ is a candidate approximant for small 
$t\in\RR$ hence the map $t\mapsto \|f-g_0-th\|^2_{L^2(I)}$ 
has a minimum at $t=0$. Differentiating under the
integral sign and equating the derivative to zero
yields $2\mbox{\rm Re}<f-g_0,h>_I=0$ whence
$<f-g_0,h>_I=0$ upon replacing $h$ by $ih$. Letting $h=e^{ik\theta}$ for
$k\in\NN$ we see that $(f-g_0)\vee 0$ 
lies in ${\bar H}_0^2$, hence it is identically zero because it vanishes on
$J$. Thus $f={g_0}_{|_I}$ as was to be shown.
\hfill \boite\\

Theorem \ref{existence} entails that the constraint 
$\|g\|_{L^\infty(J)}\leq 1$ in Problem (\ref{eq-approx}) is saturated unless
$f={g_0}_{|_I}$. If the boundary of $I$ has measure zero, 
much more in fact is true:
\begin{theorem}
Assume that $\ell(\partial I)=0$ and let $g_0$ be the solution to 
Problem {\rm (\ref{eq-approx})}.
Then $|g_0|=1$ a.e. on $J$ unless
$f=g_{|_I}$ for some $g\in H^2$ such that 
$\|g\|_{L^\infty(J)}\leq1$.
\label{thm1}
\end{theorem}
It would be interesting to know how much the assumption 
$\ell(\partial I)=0$ can be relaxed in the above statement.
Reducing Problem $BEP$ (\ref{problem}) to Problem $BEP_{2, \infty}$ 
(\ref{eq-approx}) as before, we obtain as a corollary:
\begin{corollary}
\label{corsatM}
Assume that $\ell(\partial I)=0$ and let $g_0$ be the solution to 
Problem {\rm (\ref{problem})}. If $\log M\in L^1(J)$,
then $|g_0(e^{i\theta})|=M(e^{i\theta})$ a.e. on $J$ unless
$f=g_{|_I}$ for some $g\in H^2$ such that 
$|g(e^{i\theta})|\leq M(e^{i\theta})$ a.e. on $J$.
\end{corollary}
To prove Theorem \ref{thm1}
we establish three lemmas, the second of which will be of repeated use in
the paper.
\begin{lemma}
Let $E\subset\TT$ be infinite and $K_1\subset\TT$ be a compact set such that
$\overline{E}\cap K_1=\emptyset$. If we define a collection ${\cal R}$ of 
rational functions in the variable $z$ by
\begin{equation}
\label{descS}
{\cal R}= \left\{\,c_0 + \,i\,
\sum_{k=1}^n c_k \frac{e^{i \, \psi_k} + 
  z}{e^{i \, \psi_k} - z} \, ; \right.
\end{equation}
\[
\left.
\ c_0 \, , \ c_k \in \RR \, , ~\
e^{i \, \psi_k} \in E\,, ~~1 \leq k \leq n,~~ 
n \in\NN \,
\right\} \,, 
\]
then ${\cal R}_{|_{K_1}}$ is uniformly dense in $C_\RR(K_1)$.
\label{stw}
\end{lemma}
{\sl Proof.}
It is elementary to check that members of ${\cal R}$ are real-valued a.e. 
on $\TT$. Also, it is enough to assume that 
$E$ consists of a sequence
$\{e^{i\psi_k}\}_{k\in\NN}$ that converges in $\TT$ to some
$e^{i\psi_\infty}$. 
We work over the real axis where computations are
slightly simpler, and for this we consider the
Möbius transform:
\[
\varphi(z) = i \, \frac{e^{i\psi_\infty} + z}{e^{i\psi_\infty} - z} \, ,
\]
that maps $\TT$ onto $\RR \cup \{\infty\}$ with
$\varphi(e^{i\psi_\infty})=\infty$. 
Set $K_2=\varphi(K_1)$, and note that it is compact in $\RR$
since $e^{i\psi_\infty} \notin K_1$.
Let ${\cal R}_\RR$ denote the collection of all functions 
$r\circ\varphi^{-1}$ as $r$ ranges over $\mathcal{R}$.
We are now left to prove that the restrictions 
to $K_2$ of functions in ${\cal R}_\RR$ are uniformly dense in
$C_{\RR}(K_2)$. For this,
we put $t_k=\varphi(e^{i\psi_k})$
and, denoting by $t=\varphi^{-1}(z)$ the independent variable in $\RR$,
we compute from (\ref{descS}) that
\[
{\cal R}_\RR =\{a_0 + \sum_{k=1}^n \frac{b_k}{t - t_k } \, , \ a_0 \, ,
\ b_k \in \RR\,,~~  1
\leq k \leq n  \, ,~~n \in \NN \,\}, 
\]
that is to say ${\cal R}_\RR$
is the set of real rational functions bounded at infinity,
each pole of which is simple
and coincides with some $t_k$. Thus if $P_{\RR,n}$
stands for the space of real polynomials of degree at most $n$, we get
\[
{\cal
R}_\RR = \left\{\frac{p_n(t)}{\prod_{k=1}^n (t - t_k)} \, , \ p_n \in
P_{\RR,n},~ \,\,1 \leq k \leq n,~~n\in\NN\right\} \, ,
\]
where the empty product is 1.
We claim that to each $\epsilon >0$ and
$p \in P_{\RR,n}$ there exists $r \in {\cal R}_\RR$ such that 
$$ ||r-p||_{L^{\infty}(K_2)} \leq \epsilon ,$$
and this will achieve the proof since $P_{\RR,n}$ is dense in $C_{\RR}(K_2)$
by the Stone-Weierstrass theorem. To establish the claim, let $U$ be 
a neighborhood of $0$ in $\RR^n$ such that 
$$ \forall (x_1 \dots x_n) \in U,\,\, ~\left|1-\frac{1}{\prod_{k=1}^{n}
  (1-x_k)}\right| \leq  \frac{\epsilon}{1+||p||_{L^{\infty}(K_2)}}.
$$
Next, pick $n$ distinct numbers $t_{k_1}, \dots ,t_{k_n}$
so large in modulus that $t/t_{k_j}\in U$ for $t \in K_2$ and $1\leq j\leq n$;
this is certainly possible since $K_2$ is compact whereas $|t_k|$ tends to
$\infty$ because $e^{i\psi_k}\to e^{i\psi_\infty}$. 
Finally, set 
$$ r(t)=\frac{p(t)}{\prod_{j=1}^{n} (1-\frac{t}{t_{k_j}})}.$$
Clearly $r$ belongs to ${\cal R}_\RR$, and 
$$||p-r||_{L^{\infty}(K_2)} \leq ||p||_{L^{\infty}(K_2)}\,\left\|1-
\frac{1}{\prod_{j=1}^{n} (1-\frac{t}{t_{k_j}})}\right\|_{L^{\infty}(K_2)} 
\leq \epsilon $$ 
as claimed.
\hfill \boite\\
\begin{lemma}
Let $f\in L^2(I)$ and $g_0$ be the solution to 
problem {\rm (\ref{eq-approx})}. For $h$ a 
real-valued Dini-continuous function on $\TT$ 
\emph{supported on the interior} $\stackrel{\circ}{I}$ \emph{of}
$I$, let 
\begin{equation}
b(z) = \frac{1}{2 \, \pi} \, \int_{I} \frac{e^{it} + 
z}{e^{it} - z} \, h(e^{it}) \, d t \, 
,~~~~z\in\DD,
\label{def-b}
\end{equation}
be the Riesz-Herglotz transform of $h$.
Then $b$ is continuous on $\overline{\DD}$, and moreover
\begin{equation}
\mbox{\rm Re} \,
<(f - g_0)\, \overline{g}_0 \, , \, b>_{I} = 0 \, .
\label{orthob}
\end{equation}
\label{orthopart}
\end{lemma}
{\sl Proof.}
It follows from (\ref{modDini}) that $b$ continuous on
$\overline{\DD}$. For $\lambda \in \RR$, consider the function
\[ 
\omega_{\lambda}(z) = \exp \lambda \, b(z) \, 
,~~~~z\in\DD, 
\] 
which is the outer function in $H^\infty$ whose modulus is equal to 
$\exp \lambda\, h$. Since $|\omega_\lambda|=1$ on $J$,
the function $g_0\, \omega_{\lambda}$ is a candidate approximant
in problem {\rm (\ref{eq-approx})} thus 
$\lambda \to \|f - g_0\, \omega_{\lambda}\|_{L^2(I)}^2$ reaches a
minimum at $\lambda = 0$. 
By the boundedness of $b$, we may differentiate this function with respect to
$\lambda$ under the integral sign, and equating the derivative to 0 at
$\lambda=0$ yields (\ref{orthob}).
\hfill \boite\\
\begin{lemma}
\label{supcst}
Let $f\in L^2(I)$ and $g_0$ be the solution to 
Problem {\rm (\ref{eq-approx})}. Then $(f - g_0) \, \overline{g}_0$ has
real mean on $I$:
\begin{equation}
\mbox{\rm Re} \, <(f - g_0) \, \overline{g}_0 \,,\, i>_{I} = 0.
\label{cstes}
\end{equation}
\end{lemma}
{\sl Proof.}
For each $\alpha \in [- \pi, \pi]$,  
the function $g_0 \, e^{i\alpha}$ belongs to $H^2$ and is a candidate
approximant in {\rm (\ref{eq-approx})} since it has the same modulus as $g_0$.
Hence the function $\alpha \to \|f - g_0 \, e^{i\alpha}\|_{L^2(I)}$ reaches 
a minimum at $\alpha  =0$, and differentiating under the integral
sign yields (\ref{cstes}).
\hfill \boite\\

{\sl Proof of  \rm{Theorem \ref{thm1}}.}
 Since $\partial J=\partial I$ has measure zero, it is equivalent to show that
$|g_0|=1$ a.e. on $\stackrel{\circ}{J}$.
Let \[
E = \{ e^{i \theta} \in \stackrel{\circ}{J} \, , 
\ |g_0(e^{i \theta})| < 1 \} \, ,
\]
and assume for a contradiction that $\ell(E) > 0$. 
By countable additivity, there is $\varepsilon > 0$ such that
\[
E_\varepsilon = \{ e^{i \theta} \in \stackrel{\circ}{J} \, , 
\ |g_0(e^{i \theta})| \leq
1 - \varepsilon\} 
\]
has strictly positive measure. Hence by inner regularity of
Lebesgue measure, there is a compact set $K\subset E_\varepsilon$ 
such that $\ell(K)>0$, and since $K\subset \stackrel{\circ}{J}$
it is at positive distance from $I$, say, $\eta$.
For $\lambda \in \RR$ and $F$ a
measurable subset of $K$, let $w_{\lambda,F}$
be the outer function whose modulus is $\exp \lambda$ on
$F$, and 1 on $\TT \setminus F$.
By definition $w_{\lambda,F}(z)= \exp \left\{\lambda \, A_F(z)\right\}$, where
\begin{equation}
A_F(z) = \frac{1}{2 \, \pi} \, \int_{F} \frac{e^{it} + 
z}{e^{it} - z} \, \, d t \, ,~~ ~~z\in\DD
\label{aeps-def}
\end{equation}
is the Riesz-Herglotz transform of $\chi_F$.
For $\lambda < \log (1 / (1 -\varepsilon))$ the 
function $g_0 \, w_{\lambda,F}$ belongs to $H^2$ and 
satisfies $|g_0 \,w_{\lambda,F}| \leq 1$ a.e. on $J$ so that,
by definition of $g_0$, the function
$\lambda \to \|f -
g_0\, w_{\lambda,F}\|_{L^2(I)}$ reaches a minimum at 
$\lambda = 0$. From (\ref{aeps-def}), we see that $A_F$ is uniformly 
bounded on $I$ because 
$|e^{it}-e^{i\theta}|\geq\eta>0$
whenever $e^{it}\in F$ and $e^{i\theta}\in I$.
Therefore we may differentiate under the integral sign to compute the
derivative of
$\|f - g_0\, w_{\lambda,F}\|_{L^2(I)}^2$ with respect 
to $\lambda$, which gives us
\[
-2\mbox{\rm Re} \, <f - g_0 \, \exp \{\lambda \, A_F\}\,,\,  g_0 \,
A_F \, \exp \{\lambda \, A_F\}>_{I} \, .
\]
Since the latter must vanish at $\lambda = 0$ we obtain
\begin{equation}
\mbox{\rm Re} \, <f - g_0 \,,\, g_0 \, A_F >_{I} = \mbox{\rm Re} \,
<(f - g_0)\, \overline{g}_0 \, , \, A_F >_{I} = 0 \, .
\label{aepsilon}
\end{equation}
Let $e^{it_0}$ be a density point of $K$ and $I_l$ denote
the arc centered at $e^{it_0}$ of length $l$, so that
$\ell(I_l\cap K)/l\to1$ as $l\to0$. In particular
$\ell(I_l\cap K)\neq0$ for sufficiently small $l$. Noting that
\begin{equation}
\label{comparfrac}
\left|\frac{e^{it}+e^{i\theta}}{e^{it}-e^{i\theta}}-
\frac{e^{it_0}+e^{i\theta}}{e^{it_0}-e^{i\theta}}\right|
\leq 2l/\eta^2~~~~\mathrm{for}~e^{it}\in I_l\cap K,~~e^{i\theta}
\in I,
\end{equation}
and observing that $(f-g_0)\overline{g}_0\in
L^1(I)$, we get from (\ref{aepsilon})-(\ref{comparfrac}) that
\begin{equation}
\label{annulfrac}
\mbox{\rm Re} \, <(f - g_0) \overline{g}_0\,, \, 
\frac{e^{it_0}+e^{i\theta}}{e^{it_0}-e^{i\theta}}>_{I} 
\end{equation}
\[
= 
\lim_{l\to0}\mbox{\rm Re} \,
<(f - g_0)\, \overline{g}_0 \, , \, \frac{2\pi}{\ell(I_l\cap K)}A_{I_l\cap K} >_{I} = 
0 \, .
\]
Thus, if we let $\mathcal{D}_K$ denote the set of density points of $K$, we may
capsulize (\ref{annulfrac}) and (\ref{cstes}) by saying that 
$(f - g_0) \, \overline{g}_0$ is orthogonal to the \emph{real} vector space 
\[
{\cal S}_K = \{i\,c_0 + \,
\sum_{k=1}^n c_k \frac{e^{i \, \phi_k} + 
  z}{e^{i \, \phi_k} - z} \, , \ c_0 \, , \ c_k \in \RR \, , \
e^{i \, \phi_k} \in \mathcal{D}_K \, , ~\ 1 \leq k \leq n,~~ 
n \in\NN \,\} \, 
\]
for the \emph{real} scalar product $\mbox{\rm Re}\,<~,~>_I$. Since
$\ell(\partial I)=0$ we can replace $I$ by $\bar{I}$ in
this product:
\begin{equation}
\mbox{\rm Re} \, <(f - g_0)\, \overline{g}_0 ,  \, r>_{\bar{I}} \,=\, 0 \,
, ~~~~\forall \, r \in {\cal S}_K \, .
\label{Sn}
\end{equation}
As $\ell(K)>0$ and almost every point of $K$ is a density point by 
Lebesgue's theorem \cite[sec. 7.12]{Rudin}, the set  $\mathcal{D}_K$
is certainly infinite. Moreover, since $K\subset\stackrel{\circ}{J}$,
we have that $\overline{I}\cap \overline{\mathcal{D}}_K=\emptyset$.
Now, Lemma \ref{stw} with $E=\mathcal{D}_K$ and
$K_1=\overline{I}$ implies in view of (\ref{Sn}) that
\begin{equation}
\label{orthoimc}
\mbox{\rm Re} \, <(f - g_0)\, \overline{g}_0\, ,\,  i \phi>_{\bar{I}}\,=\,0\,,~~~~\forall
\phi\in C_\RR(\bar{I}).
\end{equation}
By Riesz duality it follows that
$(f-g_0)\, \overline{g}_0$ is real-valued a.e. on 
$\bar{I}$.
In particular, if $h$ is a  Dini-continuous real function supported on 
$\stackrel{\circ}{I}$, (\ref{orthoimc}) holds with 
$\phi=\widetilde{h}_{|_{\bar{I}}}$. Hence
by Lemma \ref{orthopart}, where $I$ can be replaced by
$\bar{I}$, 
\begin{equation}
\label{orthogc}
<(f - g_0)\,
\overline{g}_0 \, , \, h>_{\bar{I}} = 0.
\end{equation}
However, by regularization, Dini-continuous --even smooth--
functions are
uniformly dense in the space of continuous 
functions with compact support on $\stackrel{\circ}{I}$
\cite[chap. 1, prop. 8]{Yosida}. Therefore
(\ref{orthogc}) in fact holds for every continuous $h$ supported on 
$\stackrel{\circ}{I}$. Consequently
$(f - g_0)\, \overline{g}_0$ must vanishes a.e. on  $\stackrel{\circ}{I}$
thus also on $I$. This implies that either $g_0= f$ a.e. on $I$ or 
$g_0 = 0$ on a set of positive measure, in which case
$g_0=0$. In any case, by Theorem \ref{existence},
$f$ is the trace on $I$ of a $H^2$-function with modulus at most 1  on $J$.
\hfill \boite\\

We now turn to the continuity of the solution to problem {\rm
  (\ref{eq-approx})} with respect to the data.
\begin{theorem}
Let $f\in L^2(I)$ and $g_0$ be the solution to 
problem {\rm (\ref{eq-approx})}.
Assume that $f^{\{n\}}$ converges to $f$ in $L^2(I)$ as $n\to\infty$,
and let $g_0^{\{n\}}$
indicate the corresponding solution to problem {\rm (\ref{eq-approx})}.
Then ${g_0^{\{n\}}}_{|_I}$ converges to ${g_0}_{|_I}$
in $L^2(I)$ and ${g_0^{\{n\}}}_{|_J}$ converges weak-* to ${g_0}_{|_J}$ in
$L^\infty(J)$. If moreover $\ell(\partial I)=0$
and $f$ is not the trace on $I$ of a $H^2$-function less than 1 in modulus a.e.
on $J$,
then $g_0^{\{n\}}$ converges to ${g_0}$ in $L^2(\TT)$.
\label{thm2}
\end{theorem}
{\sl Proof.} 
By definition $\|g_0^{\{n\}}\|_{L^\infty(J)}\leq1$, and by Theorem
\ref{existence} 
\[
\|g_0^{\{n\}}\|_{L^2(I)}\leq\|f^{\{n\}}\|_{L^2(I)},
\]
hence
$g_0^{\{n\}}$ is a bounded sequence in $H^2$. Let $g_\infty$
be a weak accumulation point and $g_0^{\{k_n\}}$ a subsequence converging
weakly to $g_\infty$ in $H^2$; {\it a fortiori} 
${g_0^{\{k_n\}}}_{|_I}$ converges weakly to ${g_\infty}_{|_I}$ in 
$L^2(I)$.
By weak (resp. weak-*) compactness of balls in $L^2(I)$
(resp. $L^\infty(J)$), we get $|g_\infty| \leq 1$ a.e. on 
$J$ and  
\[\|f-g_\infty\|_{L^2(I)}\leq\liminf_{n\to\infty}
\|f^{\{k_n\}}-g_0^{\{k_n\}}\|_{L^2(I)}.
\]
In particular $g_\infty$ is a candidate approximant, so one has
the series of inequalities:
\begin{equation}
\label{inegws}
\|f-g_0\|_{L^2(I)}\leq\|f-g_\infty\|_{L^2(I)}\leq\liminf_{n\to\infty}
\|f^{\{k_n\}}-g_0^{\{k_n\}}\|_{L^2(I)}
\end{equation}
\[
\leq\limsup_{n\to\infty}
\|f^{\{k_n\}}-g_0^{\{k_n\}}\|_{L^2(I)}.
\]
If one of these were strict, there would exist
$\varepsilon>0$ such that 
\begin{equation}
\label{petit}
\|f-g_0\|_{L^2(I)}+\varepsilon\,\leq\,
\|f^{\{k_n\}}-g_0^{\{k_n\}}\|_{L^2(I)}
\end{equation}
for infinitely many $n$. But $\|f-f^{\{k_n\}}\|_{L^2(I)}<\varepsilon/2$
for large $n$, thus for infinitely many $n$ (\ref{petit}) yields
\[
\|f^{\{k_n\}}-g_0\|_{L^2(I)}+\varepsilon/2\,\leq\,
\|f^{\{k_n\}}-g_0^{\{k_n\}}\|_{L^2(I)}
\]
contradicting the definition of $g_0^{\{k_n\}}$. Therefore equality
holds throughout in (\ref{inegws}), whence $g_\infty=g_0$ by the uniqueness part of Theorem \ref{existence}. Equality in (\ref{inegws}) is also to the effect that
\[
\lim_{n\to\infty}f^{\{k_n\}}-g_0^{\{k_n\}}=f-g_0~~~~{\rm in}~~L^2(I)
\]
because the norm of the weak limit is not less than the limit of 
the norms. Refining $k_n$ if necessary, we can assume in addition that
${g_0^{\{k_n\}}}_{|_J}$ converges weak-* to some $h$ in 
$L^\infty(J)$,
and since we already know that it converges weakly to ${g_0}_{|_J}$ in 
$L^2(J)$ we get $h={g_0}_{|_J}$.
Finally if $\ell(\partial I)=0$, we deduce from Theorem \ref{thm1}
that $|g_0|=1$ a.e. on $J$ hence
${g_0^{\{k_n\}}}_{|_J}$ converges to ${g_0}_{|_J}$ in 
$L^2(J)$ for again the norm of the weak limit is not less than the 
limit of the norms. Altogether we have shown that any sequence meeting the 
assumptions contains a subsequence satisfying the conclusions,
which is enough to prove the theorem.
\hfill \boite\\

To conclude this section, we prove that if $f$ has more summability than 
required, then so does $g_0$. 
\begin{proposition}
Assume that $f\in L^p(I)$ for some finite $p>2$.
If $g_0$ denotes the solution to problem (\ref{eq-approx})
and if $\ell(\partial I)=0$, then $g_0\in H^p$ and 
$\|g_0\|_{L^p(I)}\leq (1+K_{p/2})\|f\|_{L^p(I)}$.
\label{thm3}
\end{proposition}
{\sl Proof.}
Let $h$ be a Dini-continuous real-valued function supported in
$\stackrel{\circ}{I}$, and $b$ his Riesz-Herglotz transform.
Since $b$ has real part $h$ on $\TT$, Lemma \ref{orthopart} gives us 
\begin{equation}
\label{egalLp}
<|g_0|^2\, , \, h>_{I} \,=\, \mbox{\rm Re} \,
<f \overline{g}_0 \, , \, b>_{I}. 
\end{equation}
Using H\"older's inequality in (\ref{egalLp}) and observing that
$\|g_0\|_{L^2(I)}\leq\|f\|_{L^2(I)}
\leq\|f\|_{L^p(I)}$ in view of Theorem \ref{existence} and the fact that $p>2$
while $\ell(I)<1$, we obtain, with $1/p+1/2+1/{s_0}=1$:
\[
\left|<|g_0|^2\, , \, h>_{I} \right|\,\leq\,
\|f\|_{L^p(I)}\,\|g_0\|_{L^2(I)}\,\|b\|_{L^s(I)}\leq
\|f\|^2_{L^p(I)}\,\|b\|_{L^{s_0}(I)} \, .
\]
Thus, because the conjugation operator has norm $K_{s_0}$
on $L^{s_0}(\TT)$ while $h$ is supported on $I$, 
we get {\it a fortiori} 
\begin{equation}
\label{inegLpg}
\left|<|g_0|^2\, , \, h>_{I} \right|\,\leq\,
(1+K_{s_0})\|f\|^2_{L^p(I)}\,\|h\|_{L^{s_0}(I)}.
\end{equation}
Now, Dini-continuous functions supported on $\stackrel{\circ}{I}$ are dense in
$L^{s_0}(\stackrel{\circ}{I})$, hence also in $L^{s_0}(I)$ as $\ell(\partial I)=0$.
Therefore (\ref{inegLpg}) implies by duality
\begin{equation}
\label{majg2p}
\left\|g_0\right\|_{L^{p_1}(I)}\leq (1+K_{s_0})^{1/2}\|f\|_{L^p(I)},~~~~~~1/p_1=(1/p+1/2)/2.
\end{equation}
H\"older's inequality in (\ref{egalLp}), using this time (\ref{majg2p}) instead of 
$\|g_0\|_{L^2(I)}\leq\|f\|_{L^p(I)}$, strengthens  (\ref{inegLpg}) to
\[
\left|<|g_0|^2\, , \, h>_{I} \right|\,\leq\,
(1+K_{s_0})^{1/2}(1+K_{s_1})\|f\|^2_{L^p(I)}\,\|h\|_{L^{s_1}(I)}, \ 1/p+1/p_1+1/s_1=1,
\]
which gives us by duality 
\[
\left\|g_0\right\|_{L^{p_2}(I)}\leq (1+K_{s_0})^{1/4}\,(1+K_{s_1})^{1/2}\,\|f\|_{L^p(I)},~~~~~~
1/p_2=(1/p+1/p_1)/2.
\]
Set $1/p_k=(1/p+1/p_{k-1})/2$ and $1/p+1/p_{k}+1/s_{k}=1$.
Iterating this reasoning yields by induction
\begin{equation}
\label{prodKs}
\|g_0\|_{L^{p_k}(I)}\leq  \|f\|_{L^p(I)}
\,\,\Pi_{j=0}^{k-1}(1+K_{s_j})^{1/2^{k-j}}\,.
\end{equation}
As $k$ goes large $p_k$ increases to $p$ and
$K_{s_k}=K_{p_{k+1}/2}$ decreases to $K_{p/2}$. 
Hence the product on the right of (\ref{prodKs}) becomes arbitrarily close to
$1+K_{p/2}$, and the result now follows on letting $k\to+\infty$.
\hfill \boite\\

In problem (\ref{eq-approx}), it would be interesting to know whether
$g_0\in BMOA$ when $f\in L^\infty(I)$ and $\ell(\partial I)=0$.
\newpage
\section{The critical point equation}
\label{critpsec}
In any convex minimization problem, the solution 
is characterized by a variational inequality saying that the
{\it criterium } increases with admissible increments of the variable. 
If the problem is smooth,
infinitesimal increments span a half-space whose boundary hyperplane is tangent to the 
admissible set, and the variational inequality 
becomes an equality asserting that the derivative of the objective function
is zero on that hyperplane.
This equality, sometimes called a \emph{critical point equation},
expresses that the vector gradient of the objective function in the ambient space
lies orthogonal to the constraint; 
this vector is an implicit parameter of the critical point equation,
known as a \emph{Lagrange parameter}. 

In problem
(\ref{eq-approx}) the variational inequality is (\ref{varineg2}). However,
the non-smoothness of the $L^\infty$-norm makes it {\it a priori}
unclear whether
a critical point equation exists. It turns out that it does, at least when
$\ell(\partial I)=0$.
\begin{theorem}
Assume that $f \in L^2(I)$ is not the trace on $I$ of a $H^2$-function 
of modulus less than or equal to $1$ a.e on $J$, and suppose further that
$\ell(\partial I)=0$. Then, $g_0\in H^2$ is the solution 
to problem {\rm (\ref{eq-approx})} if, and only if, the following 
two conditions hold:
\begin{itemize}
\item[(i)] $|g_0(e^{i\theta})|=1$ for a.e. $e^{i\theta} \in J$, 
\item[(ii)] there exists
a non-negative function $\lambda \in L_{\RR}^1(J)$ such that, 
\begin{equation}
({g_0}_{|_I}-f)\,\vee\,\lambda\, {g_0}_{|_J}\,\in\,{\bar H}_0^1.
\label{carac_bepg}
\end{equation}
\end{itemize} 
Moreover, if $f\in L^p(I)$ for some $p$ such that $2<p<\infty$,
then $\lambda\in L_{\RR}^p(J)$.
\label{thm4}
\end{theorem}
{\bf Remark:} Note that (\ref{carac_bepg}) is equivalent to saying that
$({g_0}_{|_I}-f)\,\vee\,\lambda\, {g_0}_{|_J}\in L^1(\TT)$ and
\begin{equation}
\label{caracb-bepg}
{\bf P}_+\left(\,({g_0}_{|_I}-f)\,\vee\,\lambda\, {g_0}_{|_J}\,\right)\,=\,0
\end{equation}
which is the critical point equation proper, with Lagrange
parameter $\lambda$. Observe that $\log\lambda\in L_{\RR}^1(J)$,
otherwise the ${\bar H}_0^1$-function
$({g_0}_{|_I}-f)\,\vee\,(\lambda \,{g_0}_{|_J})$ would be zero
hence $f={g_0}_{|_I}$, contrary to the hypothesis.\\

To prove Theorem \ref{thm4}, we need two lemmas the first of which
stands somewhat dual to Lemma \ref{orthopart}:
\begin{lemma}
Let $f\in L^2(I)$ and $g_0$ be the solution to 
problem {\rm (\ref{eq-approx})}. If $h$ is a non-negative
function in $L^\infty(\TT)$ which is supported on
$\stackrel{\circ}{J}$, and if
\begin{equation}
a(z) = \frac{1}{2 \, \pi} \, \int_{J} \frac{e^{i\theta} + 
z}{e^{i\theta} - z} \, h(e^{i\theta}) \, d \theta \, 
,~~~~z\in\DD,
\label{def-a}
\end{equation}
denotes its Riesz-Herglotz transform, then $a$ is continuous on 
$\overline{I}$ and we have that
\begin{equation}
\mbox{\rm Re} \,
<(f - g_0)\, \overline{g}_0 \, , \, a>_{I} ~\geq~ 0 \, .
\label{orthoa}
\end{equation}
\label{orthopos}
\end{lemma}
{\sl Proof.} 
Since $h$ is supported in $\stackrel{\circ}{J}$, it is clear from the
definition that $a$ is continuous on $\overline{I}$.
For $t \in \RR$, let us put
\[ 
w_{t}(z) = \exp t \, a(z) \, 
,~~~~z\in\DD, 
\] 
which is the outer function in $H^\infty$ whose modulus is equal to 
$\exp \{t\, h\}$. As $h\geq0$,
the function $g_0\, w_{t}$ is a candidate approximant
in problem {\rm (\ref{eq-approx})} when $t\leq0$. Since
$t \to \|f - g_0\, w_{t}\|_{L^2(I)}^2$ can be differentiated with
respect to $t$ under the integral sign by the boundedness of $a$ on $I$, 
its derivative at $t=0$ must be non-positive by the minimizing property of $g_0$.
But this derivative is just $-2\mbox{\rm Re} \,
<(f - g_0)\, \overline{g}_0 \, , \, a>_{I}$.
\hfill \boite\\

Our second preparatory result is of technical nature:
\begin{lemma}
Assume that $f\in L^2(I)$ and let $g_0$ be the solution to 
problem {\rm (\ref{eq-approx})}. If $f\neq {g_0}_{|_I}$ and
$\ell(\partial I)=0$, then there exists 
a unique $\lambda\in L^1_{\RR}(J)$ such that
\begin{equation}
\label{apH2pre}
({g_0}_{|_I}-f)\, {\overline{g}_0}_{|_I}\vee \lambda\,\,\in {\bar
  H}_0^1\,.
\end{equation}
\label{prelimnv}
Necessarily $\lambda\geq0$ a.e. on $J$, and if $f\in L^\infty(I)$
then $\lambda\in L^p(J)$ for $1< p<\infty$. If
$f^{\{n\}}\in L^\infty(I)$ converges to $f$ in 
$L^2(I)$ while $g_0^{\{n\}}$ is the corresponding
solution to problem {\rm (\ref{eq-approx})},
and if we write by {\rm (\ref{apH2pre})}
\begin{equation}
\label{apH2pren}
\left({g_0}_{|_I}^{\{n\}}-f^{\{n\}}\right)\, {\overline{g}_0}_{|_I}^{\{n\}}\vee 
\lambda^{\{n\}}\,\,\in {\bar
  H}_0^1,~~~~~~~~{\rm with}~~\lambda^{\{n\}}\in L^1_{\RR}(J),
\end{equation}
then the
sequence of concatenated functions in {\rm (\ref{apH2pren})} converges 
weak-* in ${\bar   H}_0^1$ to the function
{\rm (\ref{apH2pre})}.
\end{lemma}
{\sl Proof.}
The uniqueness of $\lambda$ is clear because if $\lambda'\in L^1_{\RR}(J)$
satisfies (\ref{apH2pre}), then 
$0\vee (\lambda-\lambda')\,\in{\bar  H}_0^1$ so that $\lambda=\lambda'$.
To prove the existence of $\lambda$, assume first that
$f\in L^\infty(I)$ and fix $p\in(2,\infty)$. By 
proposition \ref{thm3} and H\"older's inequality, we know that 
$(g_0-f)\, \overline{g}_0\in L^p(I)$. For $h$ a real-valued
function in $L^q(J)$ where $1/q=1-1/p$, let $a$ be the Riesz-Herglotz 
transform of $0\vee h$ given by (\ref{def-a}) and put
\begin{equation}
\label{defL}
\mathcal{L}(h)=\mbox{\rm Re}<(f-g_0)\, \overline{g}_0 \, , \, a>_{I}.
\end{equation}
As  $0\vee h$ vanishes on $I$ by construction, it is clear that
\[
\mathcal{L}(h)=\mbox{\rm Re}<(f-g_0)\, \overline{g}_0 \, , \, \widetilde{0\vee
  h}>_{I},\]
and since the conjugation operator is bounded by $K_q$ on $L^q_{\RR}(\TT)$,
we obtain from Hölder's inequality
\[
\left|\mathcal{L}(h)\right|\,\leq\,K_q\left\|(f - g_0)\,
  \overline{g}_0\right\|_{L^p(I)} \|h\|_{L^q(J)}\,.
\]
Thus $\mathcal{L}$ is a continuous linear form on $L^q_{\RR}(J)$
and there exists $\lambda\in L^p_{\RR}(J)$
such that
\begin{equation}
\label{eqL}
\mathcal{L}(h)=<\lambda \, , \, h>_{J}, ~~~~h\in L^q(J).
\end{equation}
By Lemma \ref{orthopos}, $\mathcal{L}$ 
is a positive functional on bounded 
functions supported on
$\stackrel{\circ}{J}$. Hence $\lambda\geq0$ a.e. on $\stackrel{\circ}{J}$
thus also on $J$ since $\ell(\partial J)=\ell(\partial I)=0$.
As $\mbox{\rm Re}\,a=h$ and $\lambda$ is real-valued, equation 
(\ref{eqL}) gives us 
\begin{equation}
\label{eqLr}
\mathcal{L}(h)=\mbox{\rm Re}<\lambda \, , \, a>_{J}, ~~~~h\in L^q(J),
\end{equation}
and therefore, substracting (\ref{defL}) from (\ref{eqLr}), we get
\begin{equation}
\label{fonctL}
\mbox{\rm Re}<({g_0}_{|_I}-f)\, {\overline{g}_0}_{|_I}\vee \lambda \, \,, 
\,\, a>_{\TT}\,=\,0
\end{equation}
whenever $a$ is the Riesz-Herglotz transform of some $h\in L^q_{\RR}(J)$.

By regularization Dini-continuous
functions  are dense in continuous 
functions with compact support in $\stackrel{\circ}{I}$, 
so they are dense in $L^q(I)$ since $\ell(\partial I)=0$. Hence it follows
from Lemma \ref{orthopart} and the boundedness of the conjugation operator in
$L^q_{\RR}(\TT)$ that
\begin{equation}
\mbox{\rm Re} \,
<(g_0-f)\, \overline{g}_0 \, , \, b>_{I} = 0 \, .
\label{orthobps}
\end{equation}
whenever $b$ is the Riesz-Herglotz transform of some $\phi\in L^q_{\RR}(I)$.
As $\lambda$ is real-valued and ${\rm Re} \,b=0$ a.e. on $J$, 
we may rewrite (\ref{orthobps}) in the form
\begin{equation}
\mbox{\rm Re}<({g_0}_{|_I}-f)\, {\overline{g}_0}_{|_I}\vee \lambda \, \,, 
\,\, b>_{\TT}\,=\,0.
\label{orthobp}
\end{equation}
Now, by (\ref{RHt}), every $H^q$-function is the sum of three terms: 
a pure imaginary constant, the Riesz-Herglotz transform of $\phi\vee0$ for some $\phi\in
L^q_{\RR}(I)$, and the Riesz-Herglotz transform of $0\vee h$ for some $h\in
L^q_{\RR}(J)$. Therefore by (\ref{orthobp}), (\ref{fonctL}),
(\ref{cstes}) and the realness of $\lambda$, we obtain 
\[\mbox{\rm Re}<({g_0}_{|_I}-f)\, {\overline{g}_0}_{|_I}\vee \lambda \,\, , 
\,\, g>_{\TT}\,=\,0\,,~~~~\forall g\in H^q.
\]
Changing $g$ into $ig$ we see that the real part is superfluous
and letting $g(e^{i\theta})=e^{ik\theta}$ for
$k\in\NN$ we get
\begin{equation}
\label{apH2prel}
({g_0}_{|_I}-f)\, {\overline{g}_0}_{|_I}\vee \lambda\,\,\in {\bar
  H}_0^p\,.
\end{equation}
If $f$ is now an arbitrary function in $L^2(I)$
and $f^{\{n\}}$, $g_0^{\{n\}}$ are as indicated in the statement of the lemma,
we know from (\ref{apH2prel}), since $f^{\{n\}}\in L^\infty(I)$,
that there is a unique $\lambda^{\{n\}}$ meeting (\ref{apH2pren}).
By Theorem \ref{thm2} we have that $g_0^{\{n\}}\to g_0$ in $H^2$, 
hence by the Schwarz inequality
\begin{equation}
\label{convLinfL2}
\lim_{n\to\infty} \left\|\left(g_0^{\{n\}}-f^{\{n\}}\right)\, 
{\overline{g}_0^{\{n\}}}-
({g_0}-f)\, {\overline{g}_0}\right\|_{L^1(I)}=0.
\end{equation}
Besides, since $\lambda^{\{n\}}\geq0$
and the mean on $\TT$ of a $\bar{H}_0^1$-function is zero,
(\ref{apH2pren}) implies
\[
\left\|\lambda^{\{n\}}\right\|_{L^1(J)}=
\int_J\lambda^{\{n\}}(t)\,dt=
\int_I \left(f^{\{n\}} - {g_0^{\{n\}}}\right)\, {\overline{g}_0^{\{n\}}}(t)\,dt
\]
\[
\leq\left\|\left(g_0^{\{n\}}-f^{\{n\}}\right)\,
  {\overline{g}_0^{\{n\}}}\right\|_{L^1(I)}\,,
\]
and in view of (\ref{convLinfL2}) we deduce that
$\left\|\lambda^{\{n\}}\right\|_{L^1(J)}$ is bounded independently of $n$. 
Consequently the sequence 
\begin{equation}
\label{convfe}
\left({g_0^{\{n\}}}_{|_I}-f^{\{n\}}\right)\,
{\overline{g}_0^{\{n\}}}_{|_I}\vee \lambda^{\{n\}}
\end{equation}
has a weak-* convergent subsequence to some
$F$ in ${\bar H}_0^1$, regarding the latter as dual to
$C(\TT)/\mathcal{A}$ under the pairing $<~,~>_{\TT}$.
Checking this convergence on continuous functions 
supported on the interior of $I$, we conclude from (\ref{convLinfL2}) 
that $F_{|_{\stackrel{\circ}{I}}}=({g_0}_{|_I}-f)\,{\overline{g}_0}_{|_I}$
a.e. on $\stackrel{\circ}{I}$ thus also on $I$.
Therefore if we let $\lambda=F_{|_J}$, we meet (\ref{apH2pre}).
Checking the same convergence on positive functions 
supported on $\stackrel{\circ}{J}$, we deduce since
$\lambda^{\{n\}}\geq0$ that $F_{|_J}$ is non-negative.
 Finally,
since $F$ is determined by its trace $({g_0}_{|_I}-f)\,{\overline{g}_0}_{|_I}$
on $I$, there is a unique
weak-* accumulation point of the bounded sequence (\ref{convfe}) 
which is thus convergent.
\hfill \boite\\

{\sl Proof of  \rm{Theorem \ref{thm4}}.}
To prove sufficiency, assume that $g_0\in H^2$ satisfies
$(i)-(ii)$,
and let $u\in H^2$ be such that $\|u\|_{L^\infty(J)}\|\leq1$.
From (\ref{caracb-bepg}) we get
\[{\bf P}_+\left(\,0\,\vee\,\lambda\, {g_0}_{|_J}\,\right)
={\bf P}_+\left(\,(f-{g_0}_{|_I})\,\vee\,0\,\right)\in H^2,\]
thus applying Lemma \ref{lemme_integral} with $v=\lambda\, {g_0}_{|_J}$ and $g=u-g_0$, 
we obtain
\begin{equation}
\label{suitopt}
<\lambda {g_0}\,,\,u-g_0>_J\,=\,-
<{\bf P}_+\left(\,(f-{g_0}_{|_I})\,\vee\,0\,\right)\,,\,u-g_0>_{\TT}
\end{equation}
\[
\,=\,-<f-{g_0}\,,\,u-g_0>_{I}.
\]
Since ${\rm Re}<\lambda {g_0}\,,\,u-g_0>_J=
{\rm Re}<\lambda\,,\, u{\bar g}_0-1>_J$ is non-negative
because $\lambda\geq0$ and ${\rm Re}(u{\bar g}_0)\leq |u|\leq1$, we see from 
(\ref{suitopt}) that (\ref{varineg2}) is met.  

Proving necessity is a little harder. For this, let $g_0$ solve
problem \ref{eq-approx} and observe from Theorem \ref{thm1}
that $(i)$ holds. Thus we are left to prove $(ii)$;
in fact, we will show that the function 
$\lambda$ from Lemma \ref{prelimnv} meets (\ref{carac_bepg}).

Assume first that $f\in L^\infty(I)$. From Proposition \ref{thm3} we get in particular
$g_0\in H^4$, and by Lemma \ref{prelimnv} there is $\lambda\geq0$
in $L^2_{\RR}(J)$ such that (\ref{apH2pre}) holds with $\bar{H}_0^1$ replaced by
$\bar{H}_0^2$. Using $(i)$, we may rewrite this as
\begin{equation}
\label{apH2prefact}
\Bigl(({g_0}_{|_I}-f)\,\vee\, \lambda\,{g_0}_{|_j}\Bigr)\,
{\overline{g}_0} \,=\,F,~~~~~F\in {\bar
  H}_0^2\,.
\end{equation}
Let $g_0=jw$ be the 
inner-outer factorization of $g_0$. We will show that $F\in \bar{j} {\bar H}_0^2$,
\emph{and this will achieve the proof when} $f\in L^\infty(I)$.
Indeed, dividing (\ref{apH2prefact}) by $\bar{g}_0$ then yields
\begin{equation}
\label{apH2prefactm}
({g_0}_{|_I}-f)\,\vee\, \lambda\,{g_0}_{|_j}\,
 \,\in\, {{\bar w}}^{-1}{\bar
  H}_0^2
\end{equation}
which means that the concatenated function in (\ref{apH2prefactm}) is of 
the form:
\[
e^{-i\theta}\overline{g(e^{i\theta})/w(e^{i\theta})}
\]
for some $g\in H^2$. However, $g/w$ belongs to the Nevanlinna class $N^+$ 
by definition, and it also lies in  $L^2(\TT)$ 
because so does the function on the left-hand
side of (\ref{apH2prefactm}) (recall $|g_0|=1$ a.e. on $J$). Hence 
$g/w \in H^2$, implying that 
$e^{-i\theta}\overline{g(e^{i\theta})/w(e^{i\theta})}\in {\bar  H}_0^2\subset{\bar  H}_0^1$, as desired.

Let $j=bS_{\mu}$ where 
$b$ is the Blaschke product defined by (\ref{def-Blaschke}) and $S_{\mu}$ the
singular inner factor defined by (\ref{def-singulier}).
To prove that  $F\in \bar{j} {\bar H}_0^2$, it is enough 
by uniqueness of the
inner-outer factorization to establish separately that  
$F\in \bar{b} {\bar H}_0^2$ and $F\in \bar{S}_{\mu} {\bar H}_0^2$. 
To establish the former, it is sufficient to show that
$F\in {\bar b}_1 {\bar H}_0^2$ whenever $b_1$ is a finite Blaschke
product dividing $b$, {\it i.e.} such that $b=b_1b_2$ with $b_2$ a Blaschke 
product. Pick such a $b_1$ and put for simplicity $\gamma_0=b_2S_{\mu}w$, 
so that $g_0=b_1\gamma_0$. We can write $b_1=q/q^R$, where $q$ is 
an algebraic polynomial and $q^R=z^n\overline{q(1/\bar{z})}$ 
its reciprocal. 
We may assume that $q$ is monic and
$\mbox{deg}\,q>0$:
\[q(z)=z^n+\alpha_{n-1}z^{n-1}+\alpha_{n-2}z^{n-2}+\ldots+\alpha_0\,,
~~~~{\rm for~some~}n\in\NN\setminus\{0\}.
\]
When the set of monic polynomials of degree $n$ gets identified with 
$\CC^n$, taking as coordinates all the coefficients except the leading one,
the subset $\Omega$ of those polynomials whose roots lie in $\DD$ is open.
Now, if
$Q\in\Omega$ and $b_Q=Q/Q^R$ denotes the associated 
Blaschke product, the function $g=b_Q\gamma_0$ is a candidate 
approximant in Problem (\ref{eq-approx}) since $|g|=|g_0|$
on $\TT$, thus the map
\begin{equation}
\label{critB}
Q \to \|f - \gamma_0\, b_Q\|^2_{L^2(I)}
\end{equation}
reaches a minimum on $\Omega$ at $Q = q$. Let us write a generic $Q\in\Omega$
as 
\[Q(z)=z^n+a_{n-1}z^{n-1}+a_{n-2}z^{n-2}+\ldots+a_0.
\]
Because $b_Q(e^{i\theta})$ is a rational function in the variables $a_j$
whose denominator is locally uniformly bounded away from $0$ on $\TT$, 
we may differentiate (\ref{critB}) under the
integral sign with respect to $\mbox{\rm Re}\,a_j$, $\mbox{\rm Im}\,a_j$.
Since $q$ is a minimum point, equating these partial derivatives to zero at 
$(a_l)=(\alpha_l)$ yields 
\[
-2\mbox{\rm Re} \, <(f - g_0) \, \overline{\gamma}_0 \, , \, 
\left(x_j \, \frac{\partial b_Q(e^{i\theta})}{\partial \mbox{\rm Re}\,a_j}
+ y_j \, \frac{\partial b_Q(e^{i\theta})}{\partial \mbox{\rm Im}\,a_j}\right)_{\Big|_{{\stackrel{\displaystyle{a_l=\alpha_l}}{\!\!0\leq l\leq n-1}}}}\!\!\!\!>_{I}
= 0  \, , 
\]
for all $x_j \, , \ y_j \in \RR$ and every $j\in\{0,\ldots,n-1\}$. 
After a short computation, this gives us
\[
\mbox{\rm Re} \, <(f - g_0) \, \overline{\gamma}_0 \, \,, \, \,
\frac{z_j e^{ij\theta}}{q^R(e^{i\theta})}
-\frac{(x_j-iy_j) e^{i(n-j)\theta}q(e^{i\theta})}
{\left(q^R(e^{i\theta})\right)^2}>_{I}\,
= \,0  \, ,
\]
for all $z_j = x_j \, , + i \, y_j \in \CC$,
where the second argument in the above scalar product is a
function of $e^{i\theta}\in I$. Multiplying both arguments of 
this product by the unimodular function 
$\overline{b_1(e^{i\theta})}=q^R/q(e^{i\theta})$ does not affect
its value, thus
\begin{equation}
\label{orthogderBp}
\mbox{\rm Re} \, <(f - g_0) \, \overline{g}_0 \, \,, \, \,
\frac{z_j e^{ij\theta}}{q(e^{i\theta})}
-\frac{\bar{z}_j e^{i(n-j)\theta}}
{q^R(e^{i\theta})}>_{I}\,.
\in \RR \,,
\end{equation}
for all $z_j \in \CC$. In another connection, by the very definition of $q^R$, we have that
\[
\frac{ e^{i(n-j)\theta}}
{q^R(e^{i\theta})}=\frac{ e^{i(n-j)\theta}}
{e^{in\theta}\overline{q(e^{i\theta})}}=
\overline{\left(\frac{ e^{ij\theta}}{q(e^{i\theta})}\right)}
\]
hence the second argument of $<~,~>_I$ in (\ref{orthogderBp})
is pure imaginary on $\TT$, and since $\lambda$ is real a.e. on $J$
\begin{equation}
\label{orthogderBl}
\mbox{\rm Re} \, <\lambda \, \,, \, \,
\frac{z_j e^{ij\theta}}{q(e^{i\theta})}
-\frac{\bar{z}_j e^{i(n-j)\theta}}
{q^R(e^{i\theta})}>_{J}\,
= \,0  \, ,~~~~ \forall\, z_j \in \CC \,.
\end{equation}
Therefore, substracting (\ref{orthogderBp}) from
(\ref{orthogderBl}), we obtain from 
$(i)$ and (\ref{apH2prefact}) that
\begin{equation}
\label{orthogdera}
\mbox{\rm Re} \, <F \, \,, \, \,
\frac{z_j e^{ij\theta}}{q(e^{i\theta})}
-\frac{\bar{z}_j e^{i(n-j)\theta}}
{q^R(e^{i\theta})}>_{\TT}\,
= \,0  \, ,~~~~ \forall\, z_j \in \CC \,.
\end{equation}
The roots
of $q^R$ are reflected from those of $q$ across $\TT$, thus lie outside 
$\overline{\DD}$. Hence $e^{i(n-j)\theta}/q^R(e^{i\theta})\in H^2$,
and since $F\in{\bar  H}_0^2$ we see from (\ref{orthogdec}) that
(\ref{orthogdera}) simplifies to
\[
\mbox{\rm Re} \, <F \, \,, \, \,
\frac{z_j e^{ij\theta}}{q(e^{i\theta})}>_{\TT}\,
= \,0  \, ,~~~~ \forall\, z_j \in \CC \,.
\]
As $z_j$ is an
arbitrary complex number, the symbol ``$\mbox{\rm Re}$'' is redundant in this
equation, therefore $<F\,,\,e^{ij\theta}/q(e^{i\theta})>_{\TT}=0$ for
all $j\in\{0,\ldots,\,n-1\}$
and combining linearly these $n$ equations gives us
\begin{equation}
\label{orthogdb}
<F  \,, \,
\frac{p(e^{i\theta})}{q(e^{i\theta})}>_{\TT}\,
= \,0  \, ,~~~~  \forall p\in P_{n-1}\,,
\end{equation}
where $P_{n-1}$ is the space of algebraic polynomials of degree at most
$n-1$. Now, it is elementary that
\begin{equation}
\label{orthogbm}
{\bar b}_1 {\bar H}_0^2=
\frac{q^R}{q}\,{\bar  H}_0^2=\left(\frac{P_{n-1}}{q}\right)^\perp~~~~~~{\rm in}~~~
{\bar  H}_0^2,
\end{equation}
and consequently from (\ref{orthogdb}) and (\ref{orthogbm}),
we see that $F\in {\bar b}_1 {\bar H}_0^2$ as desired.

We turn to the proof that $F\in \bar{S}_{\mu} {\bar H}_0^2$, assuming
that $\mu$ is not the zero measure otherwise it is trivial.
We need introduce the inner divisors of $S_{\mu}$ which, by uniqueness of
the inner-outer factorization, are just the singular factors $S_{\mu_0}$
where $\mu_0$ is a positive measure on $\TT$ such that $\mu-\mu_0$ is
still positive. Pick such a $\mu_0$, and set $\beta_0=bS_{\mu-\mu_0}w$ so that
$g_0=S_{\mu_0}\beta_0$. For $a\in\DD$, consider the function
\[j_a(z)\,=\,\frac{S_{\mu_0}(z)+a}{1+{\bar a}S_{\mu_0}(z)}\,,~~~~z\in\DD.
\]
It is elementary to check that $j_a$ is inner, so that $\beta_0j_a$ is a
candidate approximant in problem (\ref{eq-approx})
because $|\beta_0j_a|=|g_0|$ a.e. on $\TT$. Therefore the map
\begin{equation}
\label{critj}
a \to \|f - \beta_0\, j_a\|^2_{L^2(I)}
\end{equation}
reaches a minimum on $\DD$ at $a=0$.
Since
\[
\frac{\partial j_a(z)}{\partial \mbox{\rm Re}\, a}\,=\,
\frac{1}{1+{\bar a}S_{\mu_0}(z)}-\frac{S_{\mu_0}(z)(S_{\mu_0}(z)+a)}
{\left(1+{\bar a}S_{\mu_0}(z)\right)^2}\,,
\]
\[
\frac{\partial j_a(z)}{\partial \mbox{\rm Im}\, a}\,=\,
\frac{i}{1+{\bar a}S_{\mu_0}(z)}+\frac{iS_{\mu_0}(z)(S_{\mu_0}(z)+a)}
{\left(1+{\bar a}S_{\mu_0}(z)\right)^2},
\]
are bounded for $z\in\TT$, locally uniformly with respect to $a\in\DD$, 
we may differentiate
(\ref{critj}) under the integral sign with respect to $\mbox{\rm Re}\, a$ and
$\mbox{\rm Im}\, a$, and equating both partial derivatives to zero at 
$a=0$ yields
\[
\mbox{\rm Re} \, <(f - g_0) \, \overline{\beta}_0 \, \,, \, \,
(x+iy)
-(x-iy) S_{\mu_0}^2>_{I}\,
= \,0  \, ,~~~~ \forall x \, , \ y
\in \RR \,.
\]
Multiplying both arguments of $<~,~>_I$ 
by the unimodular function $\overline{S}_{\mu_0}$ 
we get
\begin{equation}
\label{orthogderjp}
\mbox{\rm Re} \, <(f - g_0) \, \overline{g}_0 \, \,, \, \,
(x+iy)\,\overline{S}_{\mu_0}-
(x-iy) S_{\mu_0}>_{I}\,
= \,0  \, ,~~~~ \forall x \, , \ y
\in \RR \,.
\end{equation}
In another connection, as $(x+iy)\,\overline{S}_{\mu_0}-
(x-iy) S_{\mu_0}$ is pure imaginary on $\TT$ while $\lambda$ is real-valued,
\begin{equation}
\label{orthogderjl}
\mbox{\rm Re} \, <\lambda \, \,, \, \,
(x+iy)\,\overline{S}_{\mu_0}-
(x-iy) S_{\mu_0}>_{J}\,
= \,0  \, ,~~~~ \forall x \, , \ y
\in \RR \,.
\end{equation}
Substracting (\ref{orthogderjp}) from 
(\ref{orthogderjl}), we deduce from 
$(i)$ and (\ref{apH2prefact}) that
\[
\mbox{\rm Re} \, <F \, \,, \, \,
(x+iy)\,\overline{S}_{\mu_0}-
(x-iy) S_{\mu_0}>_{\TT}\,
= \,0  \, ,~~~~ \forall x \, , \ y
\in \RR. 
\]
Since $F\in{\bar  H}_0^2$ while 
$ S_{\mu_0}\in H^2$, this simplifies to
\[
\mbox{\rm Re} \, <F \, \,, \, \,
(x+iy)\,\overline{S}_{\mu_0})>_{\TT}\,
= \,0  \, ,~~~~ \forall x \, , \ y
\in \RR \,.
\]
But $x+iy$ is arbitrary in $\CC$, so the symbol ``$\mbox{\rm Re}$'' is 
redundant in the above equation and we obtain
\begin{equation}
\label{orthogderFj}
<F \, \,, \, \,
\,\overline{S}_{\mu_0})>_{\TT}\,=\,0.
\end{equation}
Put $F(e^{i\theta})=e^{-i\theta}\overline{g(e^{i\theta})}$ with $g\in H^2$,
and conjugate (\ref{orthogderFj}) after multiplying both 
arguments by $e^{i\theta}$:
\begin{equation}
\label{orthogderFjm}
<g\,,\,e^{-i\theta}S_{\mu_0}>_{\TT}\,=\,0.
\end{equation}
As $S_\mu$ is a nontrivial singular inner factor,
it follows from \cite[cor. 6.1.]{AhernClark} that the closed linear span
of the functions ${\bf P}_+(e^{-i\theta }S_{\mu_0})$ when $S_{\mu_0}$
ranges over all inner divisors of $S_{\mu}$ is equal to
$(S_{\mu}H^2)^\perp$ in $H^2$. Hence (\ref{orthogderFjm}) implies that
$g\in S_{\mu}H^2$, and therefore $F\in \overline{S}_{\mu}\,{\bar H}_0^2$
as announced.

Having completed the proof of necessity 
when $f\in L^\infty(I)$, we now remove this restriction.
Let $f\in L^2(I)$ and $f^{\{n\}}\in L^\infty(I)$ converge to $f$ in $L^2(I)$. 
Adding to
$f^{\{n\}}$  a small $L^2(I)$-function that goes to zero with $n$ if necessary,
we may assume that $f^{\{n\}}\notin H^2_{|_I}$.
With the notations
of Lemma  \ref{prelimnv}, let us put for simplicity
\begin{equation}
\label{defF}
F^{\{n\}}\stackrel{\Delta}{=}
\left({{g_0}_{|_I}^{\{n\}}}-f^{\{n\}}\right)\, {{\overline{g}_0}_{|_I}^{\{n\}}}
\vee \lambda^{\{n\}},~~~~
F\stackrel{\Delta}{=}\left({g_0}_{|_I}-f\right)\, {\overline{g}_0}_{|_I}\vee \lambda.
\end{equation}
By the first part of the proof, we can write
\begin{equation}
\label{jnm}
F^{\{n\}}=\bar{g}_0^{\{n\}}G^{\{n\}}\,, \mbox{ where }
G^{\{n\}}\stackrel{\Delta}{=}\left(
({g_0}_{|_I}^{\{n\}}-f^{\{n\}})\,\vee\,\lambda^{\{n\}}\, {g_0}_{|_J}^{\{n\}}\right)
\,\in\,
{\bar H}_0^1.
\end{equation}
Note that $\|G^{\{n\}}\|_{L^1(\TT)}$ is bounded since
$\|f^{\{n\}}-g_0^{\{n\}}\|_{L^2(I)}\leq \|f^{\{n\}}\|_{L^2(I)}$ (for the zero 
function is a candidate approximant) and  
$\|\lambda^{\{n\}} g_0^{\{n\}}\|_{L^1(J)}=\|\lambda^{\{n\}} \|_{L^1(J)}$ 
is bounded by Lemma \ref{prelimnv}.
Thus, extracting a subsequence if necessary, we may assume that 
$G^{\{n\}}$ converges
weak-* to some $G\in{\bar H}^1_0$,  and then $G^{\{n\}}(z)\to G(z)$
for fixed $z\in\overline{\CC}\setminus\overline{\DD}$
by (\ref{Cauchy2}). Moreover, still
from  Lemma  \ref{prelimnv}, we know that 
$F^{\{n\}}$ converges to $F$ weak-*  in ${\bar
  H}_0^1$, so we get by (\ref{Cauchy2}) again
that $F^{\{n\}}(z)\to F(z)$ for fixed
$z\in\overline{\CC}\setminus\overline{\DD}$. 
Finally Theorem \ref{thm2} entails that 
$\overline{g}_0^{\{n\}}\to \overline{g}_0$ in ${\bar H}^2$, hence using
(\ref{Cauchy2}) once more we get that
$\bar{g}_0^{\{n\}}(z)\to \bar{g}_0(z)$ for
fixed $z\in\overline{\CC}\setminus\overline{\DD}$.
Altogether, in view of (\ref{jnm}), this implies
\[F(z)=\lim_{n\to\infty} F^{\{n\}}(z)=\bar{g}_0(z)\,G(z)\,,~~~~
z\in\overline{\CC}\setminus\overline{\DD}\,,
\]
showing that $F/\bar{g}_0=G\in{\bar H}_0^1$.
By $(i)$ and the definition (\ref{defF}) of $F$,
this yields (\ref{carac_bepg}) and achieves the proof.
\hfill \boite\\

Using Theorem \ref{thm4} it is easy to characterize the solution
to problem (\ref{problem}). For this, we
write $L^1(M^2d\theta,J)$ to mean those functions $h$
on $J$ such that  $h M^2\in L^1(J)$.
\begin{corollary}
\label{genM}
Assume that $M\in L^2(J)$ is non-negative with $\log M\in L^1(J)$,
and that $f \in L^2(I)$ is not the trace on $I$ of an $H^2$-function 
of modulus less than or equal to $M$ a.e on $J$; suppose further that
$\ell(\partial I)=0$. Then, for $g_0\in H^2$ 
to be the solution to problem {\rm (\ref{problem})}, it is necessary
and sufficient that the following two properties hold:
\begin{itemize}
\item[(i)] $|g_0(e^{i\theta})|=M(e^{i\theta})$ for a.e. $e^{i\theta} \in J$, 
\item[(ii)] there exists
a non-negative measurable function $\lambda\in L^1(M^2d\theta,J)$, such that:
\begin{equation}
({g_0}_{|_I}-f)\,\vee\,\lambda\, {g_0}_{|_J}\,\in\,\bar{w}_M^{-1}{\bar H}_0^1,
\label{carac_bepgM}
\end{equation}
\end{itemize} 
where $w_M$ designates
the outer function with modulus $1$ a.e. on $I$ and modulus $M$ a.e. on $J$.
In particular if  $1/M\in L^\infty(J)$ 
(more generally if $\lambda M\in L^1(J)$), then
(\ref{carac_bepgM}) amounts to:
\begin{equation}
({g_0}_{|_I}-f)\,\vee\,\lambda\, {g_0}_{|_J}\,\in\,{\bar H}_0^1.
\label{carac_bepgMb}
\end{equation}
\end{corollary}
{\bf Remark:} We observe that, of necessity, 
$\log\lambda\in L^1(J)$.\\

{\sl Proof.}
Clearly $(i)$ is equivalent to $|g_0/w_M|=1$ a.e. on $J$, and since
$|w_M|^2=1\vee M^2$ we see on
multiplying (\ref{carac_bepgM}) by $\bar{w}_M$ that it is equivalent to
\[
\left(\frac{{g_0}_{|_I}}{w_M}-\frac{f}{w_M}\right)\,\vee\,\left(\lambda\,
M^2\right)\,\frac{{g_0}_{|_J}}{w_M}\,\,\in\,{\bar H}_0^1.
\]
The conclusion now follows from Theorem \ref{thm4} and
the reduction of problem (\ref{problem}) to 
problem (\ref{eq-approx}) given in section \ref{wposed}.
If $\lambda M\in L^1(J)$ so does $\lambda {g_0}_{|_J}$ by $(i)$,
and the function (\ref{carac_bepgM}) lies in
$e^{-\i\theta}\overline{N^+}\cap L^1(\TT)=\bar{H}_0^1$.
\hfill \boite\\

Relation (\ref{carac_bepgMb}) can be recast as a spectral equation for a 
Toeplitz operator, which
should be compared with those in \cite{ablinria,partI} that
form the basis of a constructive approach to $BEP_2$.
There, $\lambda$ is a constant and the operators involved are continuous.
In our case we consider the 
Toeplitz operator $\phi_{0\vee(\lambda-1)}$
\[\phi_{0\vee(\lambda-1)}(g)={\bf P}_+\left(0\vee (\lambda-1) g_{|_J}\right) \, ,
\]
having symbol $0\vee(\lambda-1)$, with values in $H^2$
and domain 
\[
\mathcal{D}=\{g\in H^2;~ \lambda g_{|_J}\in L^1(J),~
{\bf P}_+(0\vee \lambda g_{|_J})\in H^2\} \, .
\]
By Beurling's theorem \cite[chap. II, cor. 7.3]{Garnett}
$\phi_{0\vee(\lambda-1)}$ is densely defined, for $\mathcal{D}$ contains
$w_\rho H^2$ where $w_\rho$ is the outer function with modulus
$1\vee \min(1,1/\lambda)$. Note also that
$I+\phi_{0\vee(\lambda-1)}$ is injective, because if
${g}_{|_I}\vee\lambda {g}_{|_J}\in \bar{H}_0^2$ for some $g\in\mathcal{D}$
we may multiply it by $\bar{g}$ to obtain a $\bar{H}_0^1$-function $h$
which is real-valued on $\TT$ and thus identically zero by Poisson representation
of $\overline{h(1/\bar{z})}\in e^{i\theta}H^1$.
\begin{corollary}
\label{Topcar}
Let $M\in L^2(J)$ be non-negative and $1/M\in L^\infty(J)$. Assume $f \in L^2(I)$ 
is not the trace on $I$ of a $H^2$-function 
of modulus less than or equal to $M$ a.e on $J$; suppose further that
$\ell(\partial I)=0$. If $g_0$ is the solution to problem {\rm (\ref{problem})}
and $\lambda$ is as in (\ref{carac_bepgM}), then
\begin{equation}
\label{BEPinfsolM}
g_0=\bigl(I+\phi_{0\vee(\lambda-1)}\bigr)^{-1}\,{\bf P}_+(f\vee 0)\,.
\end{equation}
\end{corollary}
{\sl Proof.} From  (\ref{carac_bepgMb}) we see that
$\lambda {g_0}_{|_J}\in L^1(J)$ and that
\[{\bf P}_+(0\vee \lambda {g_0}_{|_J})={\bf P}_+((
  f-{g_0}_{|_I})\vee0)\,\in H^2, 
\]
hence $g_0\in\mathcal{D}$.
Using that $g_0={\bf P}_+(g_0)$, we now obtain (\ref{BEPinfsolM}) on
rewriting (\ref{carac_bepgMb}) as
\[{\bf P}_+\Bigl(
g_0\,+\,0\vee(\lambda-1){g_0}_{|_J}\,-\,f\vee0\Bigr)=0.\]
\hfill \boite\\

Further smoothness properties 
of $\lambda M^2\in L^1(J)$ follow from the next representation formula.
\begin{proposition}
\label{lambda}
Let $M\in L^2(J)$ be non-negative with $\log M\in L^1(J)$, and
assume that $f \in L^2(I)$ is not 
the trace on $I$ of an $H^2$-function of modulus less than or equal to $M$
a.e. on $J$.
Suppose also that $\ell(\partial I)=0$. If $g_0$ denotes the
solution to problem {\rm (\ref{problem})} and 
$\lambda\in L^1(M^2d\theta,J)$ is the non-negative function such that
(\ref{carac_bepgM}) holds, then $\lambda M^2$ extends across
$\stackrel{\circ}{J}$ to a holomorphic function $F$ on
$\overline{\CC}\setminus \overline{I}$ satisfying
\begin{equation}
\label{reflectionF}
F\left(1/\bar{z}\right)\,=\,\overline{F(z)},~~~~z\in
\overline{\CC}\setminus \overline{I}. 
\end{equation}
Moreover, we have the Herglotz-type representation:
\begin{equation}
\label{PoissonHL}
F(z)\,=\,
\frac{1}{2i \pi} \, \int_{I} \frac{e^{i\theta} + 
z}{e^{i\theta} - z} \, \mbox{\rm Im}\left\{f(e^{i\theta})\,
\overline{{g}_0(e^{i\theta})}\right\}\,d\theta \, 
,~~~~z\in\overline{\CC}\setminus\overline{I}.
\end{equation}
\end{proposition}
{\sl Proof.}
By $(i)$ of Corollary \ref{genM} we know that $|g_0|=M$ a.e. on $J$, 
hence multiplying (\ref{carac_bepgM})
by $\bar{g}_0$ we get
\begin{equation}
\left(|{g_0}_{|_I}|^2-f\mathop{\bar{g}_0}\nolimits_{|_I}\right)\,\vee\,
\lambda\,M^2\,\in\,e^{-i\theta}\,\overline{N}^+\cap L^1(\TT)={\bar H}_0^1.
\label{carac_bepgMm}
\end{equation}
Call $F$ the concatenated function on the left of (\ref{carac_bepgMm}), 
so that 
$H(z)=i\,\,\overline{F(1/\bar{z})}$ lies in $H^1$ and vanishes at zero 
since it has zero mean on $\TT$.
Clearly $H$ has real part ${\rm
  \mbox{Im}}f\mathop{\bar{g}_0}\nolimits_{|_I}\vee0$ on $\TT$, so the
Riesz-Herglotz representation (\ref{RHt}) yields:
\[
i\,\,\overline{F(1/\bar{z})}\,=\,
\frac{1}{2 \pi} \, \int_{I} \frac{e^{i\theta} + 
z}{e^{i\theta} - z} \, \mbox{\rm Im}\left\{f(e^{i\theta})\,
\overline{{g}_0(e^{i\theta})}\right\}\,d\theta \, 
,~~~~z\in\DD\,,
\]
and upon conjugating and changing $z$ into $1/\bar{z}$ we obtain
(\ref{PoissonHL})
for $z\in\overline{\CC}\setminus\overline{\DD}$. As the right-hand side
extends analytically to $\DD$ across $\stackrel{\circ}{J}$ by reflection,
(\ref{reflectionF}) follows.
\hfill \boite\\

The interpretation of $\lambda$ as a Lagrange parameter is justified by the
duality relation below. For convenience, we write 
$L^1_+(M^2d\theta,J)$ for the set of non-negative functions in
$L^1(M^2d\theta,J)$ whose logarithm lies in $L^1(J)$.
\begin{proposition}
\label{dual}
Assume that $M\in L^2(J)$ is non-negative with $\log M\in L^1(J)$,
and that $f \in L^2(I)$ is not the trace on $I$ of an $H^2$-function 
of modulus less than or equal to $M$ a.e on $J$. Suppose further that
$\ell(\partial I)=0$, and let $g_0\in H^2$ 
be the solution to Problem {\rm \ref{problem}} with $\lambda$ as in 
(\ref{carac_bepgM}). Then, it holds that
\begin{equation}
\|f-g_0\|^2_{L^2(I)}=\displaystyle \max_{\mu\in L^1_+(M^2d\theta,J)}
\ \min_{g\in H^2}
\ \|f-g\|^2_{L^2(I)}+\int_J\mu\,(|g|^2-M^2)\,d\theta
\label{lambda_dual}
\end{equation}
\[
= \displaystyle \min_{g\in H^2}\ \max_{\mu\in L^1_+(M^2d\theta,J)}
\ \|f-g\|^2_{L^2(I)}+\int_J\mu\,(|g|^2-M^2)\,d\theta.
\]
Moreover, the $\max\min$ and  the $\min\max$ 
are simultaneously met for $g=g_0$ and
$\mu=\lambda$. 
\end{proposition}
{\sl Proof.}
Let $A$, $B$ respectively denote the  
$\max\min$ and  the $\min\max$ in (\ref{lambda_dual}). Setting
$g=g_0$ for each $\mu$, we get
$\|f-g_0\|^2_{L^2(I)}\geq A$ from Corollary \ref{genM}-$(i)$. 
For the reverse inequality,
we fix $\mu=\lambda$ and we show that
\[\min_{g\in H^2}\ \|f-g\|^2_{L^2(I)}+\int_J\lambda\,(|g|^2-M^2)\,d\theta\]
is attained at $g_0$. Clearly, it is enough to minimize over those $g\in H^2$
such that $\lambda |g|^2\in L^1(J)$. Pick such a $g$, and 
for $t\in\RR$ let $g_t=g_0+t(g-g_0)$. The function
\[\Psi(t)=\|f-g_t\|^2_{L^2(I)}+\int_J\lambda\,(|g_t|^2-M^2)\,d\theta,\]
is convex and continuously differentiable
on $\RR$. Differentiating under 
the integral sign, we get
\[\Psi^\prime(t)=2\mbox{\rm Re}\left(<g_t-f,\,g-g_0>_I+<\lambda g_t,\,g-g_0>_J
\right),\]
and in particular
\begin{equation}
\label{der0psi}
\Psi^\prime(0)=2\mbox{\rm Re}\left(<({g_0}_{|_I}-f)\vee \lambda {g_0}_{|_J}\,,
\,g-g_0>_{\TT}\right)
\end{equation}
\[
=2\mbox{\rm Re}\left(<\bigl(({g_0}_{|_I}-f)\vee \lambda
  {g_0}_{|_J}\bigr)(\overline{g}-\overline{g_0})\,,\,1>_{\TT}\right).
\]
Now $(g_0-f)\vee \lambda g_0\in e^{-i\theta}\overline{N^+}$ by
(\ref{carac_bepgM}), and since 
$g-g_0\in H^2$ it also holds that $\overline{g}-\overline{g_0}\in\overline{N^+}$. Therefore
\[\bigl(({g_0}_{|_I}-f)\vee \lambda
  {g_0}_{|_J}\bigr)(\overline{g}-\overline{g_0})\in e^{-i\theta}\overline{N^+},\]
and since it belongs to $L^1(\TT)$ because $\lambda^{1/2}{g_0}_{|_J}$  and
$\lambda^{1/2}{g}_{|_J}$ both lie in $L^2(J)$, we deduce that it 
is also in $\bar{H}_0^1$. Consequently it has zero mean on $\TT$,
and we see from (\ref{der0psi}) that $\Psi^\prime(0)=0$, hence
$\Psi$ meets a minimum at $0$ by convexity. Expressing that 
$\|f-g_0\|^2_{L^2(I)}=\Psi(0)\leq\Psi(1)$ for each
$g\in H^2$ such that $\lambda |g|^2\in L^1(J)$
leads us to $\|f-g_0\|^2_{L^2(I)}\leq A$,
as desired. Thus we have proven the first equality in (\ref{lambda_dual})
and we also have shown it is an equality for $g=g_0$ and $\mu=\lambda$.

To establish that $\|f-g_0\|^2_{L^2(I)}= B$, observe first that
\[\max_{\mu\in L^1_+(M^2d\theta,J)}
\ \|f-g\|^2_{L^2(I)}+\int_J\mu\,(|g|^2-M^2)\,d\theta=+\infty\]
unless $|g|\leq M$ a.e. on $J$; indeed if $|g|>M$ on a set $E\subset J$ 
of strictly positive measure, we can set $\mu=\rho\chi_E+\varepsilon$ for
fixed $\varepsilon>0$ and
arbitrarily large $\rho$. 
Thus we may restrict the minimization in the second line of
(\ref{lambda_dual}) to those $g$ such that $|g|\leq M$ a.e. on $J$.
For such $g$ the maximum is attained when $\mu=0$, and by
definition $g_0$ minimizes $\|f-g_0\|^2_{L^2(I)}$ among them.
Moreover, when $g=g_0$, we see from Corollary
\ref{genM}-$(i)$ that $\mu$ is irrelevant in the criterion and can be 
chosen to be $\lambda$. This achieves the proof.
\hfill \boite\\

Note that Proposition \ref{dual} would still hold if we dropped the log-integrability requirement in
the definition of $L^1(M^2d\theta,J)$, for the latter was never needed in the proof. However, 
this requirement conveniently restricts the maximization space in (\ref{lambda_dual}) to a class
of $\mu$ for which one can form the outer function $w_\mu$, and this will be of use in what follows.
\newpage
\section{The dual functional and Carleman's formulas}
\label{Carlsec}
For $M\in L^2(J)$ a non-negative function such that
$\log M\in L^1(J)$ and  $f\in L^2(I)$ which is  not 
the trace on $I$ of a $H^2$-function of modulus less than or equal to $M$
a.e. on $J$, we denote by $\Phi_M$ the \emph{dual functional} in problem (\ref{problem})
acting on $L^1_+(M^2d\theta,J)$ as follows (compare \cite[sec. 4.3]{BorLew}):
\begin{equation}
\label{defcritL}
\Phi_M(\mu)=\min_{g\in H^2}
\ \|f-g\|^2_{L^2(I)}+\int_J\mu\,(|g|^2-M^2)\,d\theta,~~~~~~~~
\mu\in L^1_+(M^2d\theta,J).
\end{equation}
As an {\it infimum} of affine functions,
$\Phi_M$ is concave and upper semi-continuous with respect to $\mu$.
In view of (\ref{lambda_dual}), solving problem (\ref{problem}) amounts
to maximize $\Phi_M$ over the convex set 
$L^1_+(M^2d\theta,J)$. 
As we shall see momentarily ({\it cf.} Proposition \ref{mingdelambda}), 
the true nature of Carleman-type formulas in this
context is that they solve for the optimal $g$ in (\ref{defcritL})
whenever the {\it min} is attained. We begin with a theorem showing how Carleman's formula
solves for $g_0$ in (\ref{carac_bepgM}) as a function of $f$ and $\lambda$.
\begin{theorem}
\label{expg0}
Let $M\in L^2(J)$ be non-negative 
with $\log M\in L^1(J)$, and
assume that $f \in L^2(I)$ is not 
the trace on $I$ of a $H^2$-function of modulus less than or equal to $M$
a.e. on $J$.
Suppose that $\ell(\partial I)=0$, and let $g_0$ be the
solution to problem {\rm (\ref{problem})} while
$\lambda\in L^1(M^2d\theta,J)$ denotes the non-negative function such that
(\ref{carac_bepgM}) holds. Write
$w_{\lambda^{1/2}}$ for the outer
function with modulus $\lambda^{1/2}$ a.e. on $J$ and 
modulus $1$ a.e. on $I$. Then
\begin{equation}
\label{CarlemanM}
g_0(z)\,=\,
\frac{1}{2i\pi \,w_{\lambda^{1/2}}(z)} \, \int_{I}
\frac{w_{\lambda^{1/2}}(\xi)\,f(\xi)}
{\xi-z}\,d\xi,~~~~z\in{\DD}.
\end{equation}
Conversely, if $\lambda$ is a positive function on $J$ such that 
$\log \lambda\in L^1(J)$ and if $g_0$ defined by 
(\ref{CarlemanM}) lies in $H^2$, then $g_0$ is the solution to problem
(\ref{problem}) where $M=|g_{0|_J}|$. In this case $\lambda$ is
the function appearing in (\ref{carac_bepgM}).
\end{theorem}
{\sl Proof.} Assume $g_0$ is the solution to problem (\ref{problem}) so
that $(i)$ and $(ii)$ of Corollary \ref{genM} hold.
Dividing (\ref{carac_bepgM}) by $\bar{w}_{\lambda^{1/2}}$ and using that
$|w_{\lambda^{1/2}}|^2=1\vee\lambda$, we deduce 
\[
w_{\lambda^{1/2}}(g_0-(f\vee0))\,\,\in\,
\bar{w}_{\lambda^{1/2}}^{-1}\,\bar{w}_{M}^{-1}\,{\bar H}_0^1.
\]
Since $\lambda\in L^1(M^2d\theta,J)$, the left-hand side lies in $L^2(\TT)$
and therefore it belongs to $\bar{H}_0^2$ because the right-hand side is in
$e^{-i\theta}\overline{N^+}$ by construction. In particular
\begin{equation}
\label{etapC}
{\bf P}_+\left(w_{\lambda^{1/2}}(g_0-(f\vee0)) \right)=0.
\end{equation}
But $w_{\lambda^{1/2}}g_0\in H^2$ because it 
clearly belongs to $N^+\cap L^2(\TT)$, so that (\ref{etapC}) implies
\[
w_{\lambda^{1/2}}g_0={\bf P}_+(w_{\lambda^{1/2}}g_0)={\bf P}_+\left(
w_{\lambda^{1/2}}(f\vee0)\right).
\]
Now  (\ref{CarlemanM}) follows from this and (\ref{Cauchyproj}).
Conversely, assume that $g_0$ defined by 
(\ref{CarlemanM}) lies in $H^2$ and set $M=|g_0|_{|_J}$. 
Since $fw_{\lambda^{1/2}}\in L^2(I)$, we see from (\ref{CarlemanM})
and (\ref{Cauchyproj}) that
$g_0w_{\lambda^{1/2}}\in H^2$ and that
\[g_0w_{\lambda^{1/2}}={\bf P}_+\left(fw_{\lambda^{1/2}}\vee0\right)\]
which implies (\ref{etapC}). Thus $w_{\lambda^{1/2}}(g_0-(f\vee0))\in \bar{H}_0^2$ and
multiplying by $\bar{w}_{M}\bar{w}_{\lambda^{1/2}}\in \bar{H}^2$
yields 
\[\bar{w}_{M}\,|w_{\lambda^{1/2}}|^2\,(g_0-(f\vee0))=
\bar{w}_{M}\bigl(({g_0}_{|_I}-f)\vee\lambda {g_0}_{|_J}\bigr)\in\bar{H}_0^1\]
from which (\ref{carac_bepgM}) follows. As $(i)$ of Corollary \ref{genM}
is met by definition, $g_0$ indeed
solves for problem (\ref{problem}).
\hfill \boite\\

Theorem \ref{expg0} can be used as follows to compute the function $\Phi_M(\mu)$
introduced in (\ref{defcritL}).
\begin{proposition}
\label{mingdelambda}
Let $M\in L^2(J)$ be non-negative 
with $\log M\in L^1(J)$, and
assume that $f \in L^2(I)$ is not 
the trace on $I$ of an $H^2$-function of modulus less than or equal to $M$
a.e. on $J$.
Suppose further that $\ell(\partial I)=0$ and
let $\mu\in L^1_+(M^2d\theta,J)$. Then,
the function $\Phi_M(\mu)$ defined by (\ref{defcritL}) can be expressed as 
\begin{equation}
\label{expcritmin}
\Phi_M(\mu)=\left\|{\bf P}_{-}\left(fw_{\mu^{1/2}}\vee0\right)
\right\|^2_{L^2(\TT)}-\left\|\mu^{1/2} M\right\|_{L^2(J)}^2.
\end{equation}
Moreover, if we set 
\begin{equation}
\label{defgmu}
g_{\mu}(z)=
\frac{1}{2i\pi \,w_{\mu^{1/2}}(z)} \, \int_{I}
\frac{w_{\mu^{1/2}}(\xi)\,f(\xi)}
{\xi-z}\,d\xi,~~~~z\in{\DD},
\end{equation}
then the infimum in the right-hand side of (\ref{defcritL}) is attained at $g=g_\mu$
whenever the latter belongs to $H^2$. In particular, this is the case when
$1/\mu\in L^\infty(J)$.
\end{proposition}
{\sl Proof.}
Assume first that $\mu$ is such that $g_\mu\in H^2$; this holds in particular when 
$1/\mu\in L^\infty(J)$, because then
$1/w_{\mu^{1/2}}\in H^\infty$ while (\ref{Cauchyproj}) shows
that the integral in (\ref{defgmu}) lies in $H^2$. From
Theorem \ref{expg0} it follows that $g_\mu$ is the solution to problem (\ref{problem})
where $M$ gets replaced by $|g_\mu|$, and $\mu$ plays the role of
$\lambda$ in (\ref{carac_bepgM}). Hence Proposition \ref{dual} implies
that $g_\mu$ is an infimizer of
$$\min_{g\in H^2}
\ \|f-g\|^2_{L^2(I)}+\int_J\mu\,(|g|^2-|g_\mu|^2)\,d\theta,$$
and since $\mu$ is kept fixed $g_\mu$ is clearly also an infimizer
of
$$\min_{g\in H^2}
\ \|f-g\|^2_{L^2(I)}+\int_J\mu\,(|g|^2-M^2)\,d\theta$$
which is just the right-hand side of (\ref{defcritL}). This proves the
second assertion of the proposition.

By (\ref{defgmu})) and (\ref{Cauchyproj}), taking into account that
$|w_{\mu^{1/2}}|=1\vee\mu^{1/2}$, what precedes 
can be reformulated as 
\[
\Phi_M(\mu)=
\displaystyle \|f-g_\mu\|^2_{L^2(I)}+\int_J\mu\,(|g_\mu|^2-M^2)\,d\theta
\]
\[= 
\displaystyle \left\|(w_{\mu^{1/2}}f\vee 0)-w_{\mu^{1/2}}g_\mu\right\|^2_{L^2(\TT)}-
\int_J\mu\,M^2\,d\theta\]
\[ 
=
\displaystyle \left\|P_{\bar{H}_0^2}\left(fw_{\mu^{1/2}}\vee0\right)
\right\|^2_{L^2(\TT)}-\left\|\mu^{1/2} M\right\|_{L^2(J)}^2.
\]
This proves (\ref{expcritmin}) when $g_\mu\in H^2$. To get it in
general we apply what we just did to the sequence $\mu_n=\mu+1/n$, observing
that $g_{\mu_n}\in H^2$ because $1/\mu_n\in L^\infty(J)$.
By monotone convergence we obtain
\begin{equation}
\label{convmoncal}
\lim_{n\to\infty}\left\|\mu_n^{1/2} M-\mu^{1/2} M\right\|_{L^2(J)}=0.
\end{equation}
Moreover, as $\log\mu_n$ decreases to $\log \mu$, we certainly have
on putting $\log^-(x)=\max\{-\log x,0\}$ and $\log^+(x)=\max\{\log x,0\}$ that
\[
\begin{array}{lll}
\log^-\mu_n&\leq&\log^-\mu\leq|\log\mu|\in L^1(J),\\
\log^+\mu_n&\leq&\log^+\left(\mu_n M^2\right)+\left|\log M^2\right|
\leq\left|\mu_nM^2-1\right|+2|\log M|\\
&\leq&(\mu+1)M^2+1+2|\log M|\in L^1(J),
\end{array}
\]
and therefore, by dominated convergence as applied to
$\log\mu_n=\log^+\mu_n-\log^-\mu_n$, we obtain
\[\lim_{n\to\infty}
\exp \, \left\{ 
\frac{1}{4 \, \pi} \, \int_{J} \frac{e^{it} + 
z}{e^{it} - z} \, \log\mu_n d t \right\}
\,=\,
\exp \, \left\{ 
\frac{1}{4 \, \pi} \, \int_{J} \frac{e^{it} + 
z}{e^{it} - z} \, \log\mu d t \right\},~~~~
z\in \stackrel{\circ}{I},
\]
in other words $w_{\mu_n^{1/2}}$ converges pointwise to 
$w_{\mu^{1/2}}$ on $\stackrel{\circ}{I}$ and therefore almost everywhere on
$I$ since $\ell(\partial I)=0$. Thus, 
appealing to dominated convergence once more, 
we get
\begin{equation}
\label{convdomcal}
\lim_{n\to\infty}\left\|fw_{\mu_n^{1/2}} -fw_{\mu^{1/2}}\right\|_{L^2(I)}=0,
\end{equation}
and from  (\ref{convmoncal}), (\ref{convdomcal}), and (\ref{expcritmin}) which 
is known to hold with $\mu$ replaced by $\mu_n$, we see that
\begin{equation}
\label{critlim1p}
\lim_{n\to\infty} \Phi_M(\mu_n)=
\left\|P_{\bar{H}_0^2}\left(fw_{\mu^{1/2}}\vee0\right)
\right\|^2_{L^2(\TT)}-\left\|\mu^{1/2} M\right\|_{L^2(J)}^2.
\end{equation}
In another connection, it is plain that 
\begin{equation}
\label{critlim2p}
\limsup_{n\to\infty} \Phi_M(\mu_n)\leq \Phi_M(\mu)\leq
\liminf_{n\to\infty} \Phi_M(\mu_n),
\end{equation}
where the first inequality comes from (\ref{convmoncal}) and the upper
semi-continuity of $\Phi_M$ in $L^1_+(M^2d\theta,J)$
while the second inequality is obvious from (\ref{defcritL}),
(\ref{convmoncal}), and the fact that $\mu\leq\mu_n$. 
Now (\ref{expcritmin}) follows from (\ref{critlim1p}) and (\ref{critlim2p}).
\hfill \boite\\

Being concave on the convex set  $L^1_+(M^2d\theta,J)$, the functional
$\Phi_M$ has a directional derivative at every point in each admissible
direction. Here, a direction $h$ is said to be admissible at $\mu\in  L^1_+(M^2d\theta,J)$
if $\mu+th\in  L^1_+(M^2d\theta,J)$ as soon as $t\geq0$ is small enough.
From a constructive viewpoint, computing this derivative is important
when designing ascent algorithms to maximize $\Phi_M$ and thus numerically solve for problem 
(\ref{problem}). The next proposition does it, under mild assumptions on $f$,
in those directions $h$ such that
$h/\mu\in L^\infty(J)$. Note since $\mu\neq0$ a.e.(for $\log\mu\in L^1(J)$) that such directions are dense
in the set of all admissible directions, hence this result  allows one indeed 
to find a direction of ascent for $\Phi_M$.
\begin{proposition}
\label{cduGateau}
Assumptions and notations being as in Proposition \ref{mingdelambda}, suppose in addition that
$|f|^2$ lies in the Zygmund class $L\log^+ L$. Let further
$h$ be a real function on $J$ such that $\|h/\mu\|_{L^\infty(J)}<1$.
Then $\mu+h\in L^1_+(M^2d\theta,J)$ and $h\in L^1(M^2d\theta,J)$. Moreover, 
defining $g_\mu$ as in (\ref{defgmu}), it holds that $h|g_\mu|^2\in L^1(J)$ and that
\begin{equation}
\label{derGat}
\left|\Phi_M(\mu+h)-\Phi_M(\mu)-
\int_J\,h(|g_\mu|^2-M^2)\,d\theta
\right|={\bf o}\left(\|h/\mu\|_{L^\infty(J)}\right),
\end{equation}
where the function ${\bf o}$, which  depends on $f$ and $\mu$ only, 
is a little $o$ of its argument near $0$. 
\end{proposition}
{\sl Proof.}
Clearly $\mu+h=\mu(1+h/\mu)\in L^1_+(M^2d\theta,J)$ whenever $\|h/\mu\|_{L^\infty(J)}<1$, 
which in turn entails $h\in L^1(M^2d\theta,J)$. In another connection, we see from
(\ref{Cauchyproj}) that (\ref{defgmu}) can be rewritten as 
\begin{equation}
\label{gmuproj}
w_{\mu^{1/2}}g_\mu={\bf P}_+(fw_{\mu^{1/2}}\vee0),
\end{equation}
and since $|w_{\mu^{1/2}}|^2=1\vee\mu$ we see that
$w_{\mu^{1/2}}g_\mu\in H^2$ hence
\[h|{g_\mu}{}_{|_J}|^2=(h/\mu)\mu|{g_\mu}{}_{|_J}|^2\in L^1(J)\]
for $h/\mu\in{L^\infty(J)}$.
Thus the integral in the left-hand side of (\ref{derGat}) is well-defined.
Next, multiplying the $\bar{H}_0^2$-function
$w_{\mu^{1/2}}g_\mu-(w_{\mu^{1/2}}f\vee0)$
by the $\bar{H}^2$-function $\overline{w_{\mu^{1/2}}g_\mu}$ yields
\[(|{g_\mu}{}_{|_I}|^2-f{\bar{g}_\mu}{}_{|_I})\vee\mu|g_\mu{}_{|_J}|^2\in \bar{H}_0^1.\]
Therefore the conjugate 
function of $(|{g_\mu}_{|_I}|^2-\mbox{\rm Re}
(\bar{f}{g_{\mu}}{}_{|_I}))\vee\mu|{g_\mu}_{|_J}|^2$
lies in $L^1(\TT)$, and by Zygmund's theorem so does the conjugate function of $|f|^2\vee0$ 
since the latter lies in $L\log^+L$ by assumption. Adding up yields
\[\widetilde{\overbrace{
\left(\frac{|g_\mu|^2_{|_I}+|f|^2}{2}+
\frac{|{g_\mu}{}_{|_I}-f|^2}{2}\right)\vee\mu|g_\mu|_{|_J}^2}}\in L^1(\TT),\]
 and since the function under brace is positive it lies in $L\log^+ L$ by the M. Riesz theorem.
{\it A fortiori} then,
\begin{equation}
\label{Rieszlog}
\left|{\bf P}_-(f{w_{\mu^{1/2}}}{}_{|_I}\vee0)\right|^2
=
\left|(f{w_{\mu^{1/2}}}{}_{|_I}\vee0)-w_{\mu^{1/2}}g_\mu\right|^2
\end{equation}
\[
=
|{g_\mu}{}_{|_I}-f|^2\vee\mu|g_\mu|_{|_J}^2
\in L\log^+L.
\]
Now, let us write
\[w_{(\mu+h)^{1/2}}(z)=w_{\mu^{1/2}}(z)
\exp \, \left\{ 
\frac{1}{4 \, \pi} \, \int_{J} \frac{e^{i\theta} + 
z}{e^{i\theta} - z} \, \log(1+h/\mu)(e^{i\theta})\, d \theta \right\}
\]
\[
=w_{\mu^{1/2}}(z)\,e^{\Delta_h(z)},
\]
where we have put for simplicity
\begin{equation}
\label{defDeltah}
\Delta_{h}(z)=
\frac{1}{4 \, \pi} \, \int_{J} \frac{e^{i\theta} + 
z}{e^{i\theta} - z} \, \log(1+h/\mu)(e^{i\theta})\, d \theta,~~~~z\in\DD.
\end{equation}
Note that $\Delta_h\in BMOA$ since $\log(1+h/\mu)\in L^\infty(J)$.
With this notation, it is straightforward that
\begin{equation}
\label{diffprinc}
\left\|{\bf P}_{-}\left(fw_{(\mu+h)^{1/2}}\vee0\right)
\right\|^2_{L^2(\TT)}-
\left\|{\bf P}_{-}\left(fw_{\mu^{1/2}}\vee0\right)
\right\|^2_{L^2(\TT)}
\end{equation}
\[
=
\left\|{\bf P}_{-}\left(fw_{\mu^{1/2}}(e^{\Delta_{h}}-1)\vee0\right)
\right\|^2_{L^2(\TT)}
\]
\[
+2{\rm Re}<{\bf P}_{-}\left(fw_{\mu^{1/2}}\vee0\right),
{\bf P}_{-}\left(fw_{\mu^{1/2}}(e^{\Delta_{h}}-1)\vee0\right)>_{\TT},
\]
and our next goal is to prove that
\[
\left|2{\rm Re}<{\bf P}_{-}\left(fw_{\mu^{1/2}}\vee0\right)~,~
{\bf P}_{-}\left(fw_{\mu^{1/2}}(e^{\Delta_{h}}-1)\vee0\right)>_{\TT}-
\int_J\,h|g_\mu|^2\,d\theta
\right|
\]
\begin{equation}
\label{dermorceau1}={\bf o}\left(\|h/\mu\|_{L^\infty(J)}\right).
\end{equation}
For this, since ${\bf P}_++{\bf P}_- ={\rm id}$, we first observe from (\ref{orthogdec}) that
\[
<{\bf P}_{-}\left(fw_{\mu^{1/2}}\vee0\right),\,
{\bf P}_{-}\left(fw_{\mu^{1/2}}(e^{\Delta_{h}}-1)\vee0\right)>_{\TT}
\]
\[
=
<{\bf P}_{-}\left(fw_{\mu^{1/2}}\vee0\right),\,
(e^{\Delta_{h}}-1)(fw_{\mu^{1/2}}\vee0)>_{\TT}
\]
\[
=<{\bf P}_{-}\left(fw_{\mu^{1/2}}\vee0\right),\,
(e^{\Delta_{h}}-1){\bf P}_{-}\left(fw_{\mu^{1/2}}\vee0\right)>_{\TT}
\]
\[=
<\left|{\bf P}_{-}\left(fw_{\mu^{1/2}}\vee0\right)\right|^2,\,
e^{\Delta_{h}}-1>_{\TT}
\]
where we used in the second equality that 
$(e^{\Delta_{h}}-1){\bf P}_{+}\left(fw_{\mu^{1/2}}\vee0\right)\in H^2$ for
$e^{\Delta_{h}}-1\in H^\infty$. Besides, 
\[{\bf P}_{-}\left(fw_{\mu^{1/2}}\vee0\right)+
{\bf P}_{+}\left(fw_{\mu^{1/2}}\vee0\right)=0~~~~a.e.~{\rm on}~J\]
which implies in view of (\ref{gmuproj}) that
\[\int_J\,h|g_\mu|^2=
\int_J \frac{h}{\mu}\,\left|{\bf P}_+\left(fw_{\mu^{1/2}}\vee0\right)\right|^2=
<\left|{\bf P}_-\left(fw_{\mu^{1/2}}\vee0\right)\right|^2, \,0\vee h/\mu>_{\TT}.\]
Altogether, the expression inside absolute values on the left-hand side of (\ref{dermorceau1}) is thus
\[<\left|{\bf P}_-\left(fw_{\mu^{1/2}}\vee0\right)\right|^2,\, {\rm Re}\left(2(e^{\Delta_{h}}-1)-
(0\vee h/\mu)\right)>_{\TT}.\]
Now, if we remark from (\ref{defDeltah}) that, on $\TT$,
we have $2\Delta_h=0\vee\log(1+h/\mu)+i\varphi$ where $\varphi$
denotes the conjugate function of $0\vee\log(1+h/\mu)$, the above quantity becomes $Q_1+Q_2$ with
\[
\begin{array}{lll}
Q_1&\stackrel{\Delta}{=}&
2\,<\left|{\bf P}_-\left(fw_{\mu^{1/2}}\vee0\right)\right|^2~,~ (\cos(\varphi/2)-1)\left(1\vee(1+h/\mu)^{1/2}
\right)>_{\TT},\\
Q_2&\stackrel{\Delta}{=}&
2\,<\left|{\bf P}_-\left(fw_{\mu^{1/2}}\vee0\right)\right|^2~,~
(1+h/\mu)^{1/2}-1-h/(2\mu)>_{J}.
\end{array}
\]
We prove separately that $Q_1$ and $Q_2$ are both ${\bf o}\left(\|h/\mu\|_{L^\infty(J)}\right)$;
here and thereafter, we use the same symbol ${\bf o}$ for different functions as this causes no confusion.
On the one hand, since there is an absolute constant $C$ such that
$\left|(1+h/\mu)^{1/2}-1-h/(2\mu)\right|<C\|h/\mu\|_{L^\infty(J)}^2$ for
$\|h/\mu\|_{L^\infty(J)}<1$, we have that
\begin{equation}
\label{majQ2}
|Q_2|\,\leq\,
2C \|f\|_{L^2(\TT)}^2\|h/\mu\|_{L^\infty(J)}^2
\end{equation}
which is indeed ${\bf o}\left(\|h/\mu\|_{L^\infty(J)}\right)$, where ``$\bf o$'' is 
independent of $\mu$. On the other hand, 
as $\cos(\varphi/2)-1\leq0$, it holds for $\|h/\mu\|_{L^\infty(J)}<1$ that
\begin{equation}
|Q_1|\,\leq\,
2\sqrt{2}<\left|{\bf P}_-\left(fw_{\mu^{1/2}}\vee0\right)\right|^2~,~ 1-\cos(\varphi/2)>_{\TT}.
\label{morceauder1}
\end{equation}
Put for simplicity 
\[
B_h\stackrel{\Delta}{=}(1-\cos(\varphi/2))/\varphi \mbox{ and } u_n\stackrel{\Delta}{=}
\min\left\{\left|{\bf P}_-\left(fw_{\mu^{1/2}}\vee0\right)\right|^2~,~n\right\},
\]
so that $B_h$ is uniformly bounded (independently of $h$)
and so is $u_n$ for fixed $n$. By monotone convergence we get
\begin{equation}
\label{convun}
<\left|{\bf P}_-\left(fw_{\mu^{1/2}}\vee0\right)\right|^2~,~ 1-\cos(\varphi/2)>_{\TT}
\end{equation}
\[
=
\lim_{n\to\infty}<u_n~,~ 1-\cos(\varphi/2)>_{\TT}
=\lim_{n\to\infty}<B_hu_n~,~ \varphi>_{\TT}.
\]
Being the product of two functions in $BMOA$,
the function
\[
(B_hu_n+i\widetilde{B_hu_n})((0\vee\log(1+h/\mu))+i\varphi)\]
certainly lies in $H^1$
and since it is real at $0$ we deduce
from the Cauchy formula  and (\ref{LlogLconj}) that
\begin{equation}
\label{majCauchy}
\left|<B_hu_n~,~ \varphi>_{\TT}\right|=\left|<\widetilde{B_hu_n}~,~\log(1+h/\mu)>_J\right|
\end{equation}
\[
\leq \|B_hu_n\|_{L\log^+L}\|\log(1+h/\mu)\|_{L^\infty(J)}.
\]
As
$|B_hu_n|\leq \left|B_h\left({\bf P}_-\left(fw_{\mu^{1/2}}\vee0\right)\right)^2\right|$ the 
same is true of their decreasing rearrangements, thus by (\ref{defLlogL}) and the inequality
$|\log(1+h/\mu)|\leq 2|h/\mu|$ which is valid for $|h/\mu|\leq1/2$, we obtain from 
(\ref{convun})-(\ref{majCauchy}) that
\[
<\left|{\bf P}_-\left(fw_{\mu^{1/2}}\vee0\right)\right|^2~,~ 1-\cos(\varphi/2)>_{\TT}
\]
\[
\leq~ 2\left\|B_h\left({\bf P}_-\left(fw_{\mu^{1/2}}\vee0\right)\right)^2\right\|_{L\log^+L}
\|h/\mu\|_{L^\infty(J)}
\] 
as soon as $\|h/\mu\|_{L^\infty(J)}<1/2$. Therefore, to prove that (\ref{morceauder1}) is 
${\bf o}(\|h/\mu\|_{L^\infty(J)})$, it is enough to show that
\begin{equation}
\label{petito}
\lim_{\|h/\mu\|_{L^\infty(J)}\to0}
\left\|B_h\left({\bf P}_-\left(fw_{\mu^{1/2}}\vee0\right)\right)^2\right\|_{L\log^+L}=0.
\end{equation}
It is easily checked that  $\|B_h\|_{L^\infty(\TT)}\leq8$, say, which entails 
for the decreasing rearrangements the inequality
\begin{equation}
\label{domR}
\left(B_h\left({\bf P}_-\left(fw_{\mu^{1/2}}\vee0\right)\right)^2\right)^*
~\leq~8\left(\left({\bf P}_-\left(fw_{\mu^{1/2}}\vee0\right)\right)^2\right)^*.
\end{equation}
Since the right-hand side of (\ref{domR}) is independent of $h$ and lies in $L\log^+L$ by
(\ref{Rieszlog}), we deduce from 
the definition (\ref{defLlogL}) of the $L\log^+L$-norm that
to each $\varepsilon>0$ there is $\eta>0$ such that
\[E\subset[0,1]~{\rm and}~\ell(E)<\eta\Longrightarrow~~\int_E
\left(B_h\left({\bf P}_-\left(fw_{\mu^{1/2}}\vee0\right)\right)^2\right)^*(t)~\log(1/t)\,dt<~\varepsilon.\]
Thus, (\ref{petito}) will hold if only
$\left(B_h\left({\bf P}_-\left(fw_{\mu^{1/2}}\vee0\right)\right)^2\right)^*$ tends to $0$ 
\emph{in measure} as $\|h_n/\mu\|_{L^\infty(J)}\to0$. 
Since a function and its decreasing rearrangement 
have the same distribution function, this is equivalent to 
\begin{equation}
\label{convmesure}
\lim_{\|h/\mu\|_{L^\infty(J)}\to0}
B_h\left({\bf P}_-\left(fw_{\mu^{1/2}}\vee0\right)\right)^2=0~~~~\mathrm{in~measure~on~}\TT.
\end{equation}
Because $|w_{\mu^{1/2}}|=1$ on $I$ and ${\bf P}_-$ is a contraction in $L^2(\TT)$,
we have the Kolmogorov estimate
\[\ell\{\xi\in\TT;~\left|{\bf P}_-\left(fw_{\mu^{1/2}}\vee0\right)\right|^2>x\}
\leq\frac{\|f\|^2_{L^2(I)}}{x},\]
hence (\ref{convmesure}) will  hold if $B_h$ alone converges to $0$ in measure.
As $|B_h|\leq C'|\varphi|$ for some absolute constant $C'$ it is sufficient to
establish that $\varphi$ in turn converges to $0$ in measure on $\TT$.
But this follows from the fact that $\varphi$ tends to $0$ in $L^p(\TT)$ 
when $\|h/\mu\|_{L^\infty(J)}\to0$ for $1<p<\infty$,
since by the M. Riesz theorem
\[\|\varphi\|_{L^p(\TT)}\leq K_p \|\log(1+h/\mu)\|_{L^p(\TT)}
\leq 2K_p \|h/\mu\|_{L^\infty(\TT)}\]
as soon as $\|h/\mu\|_{L^\infty(J)}<1/2$. This completes the proof of (\ref{dermorceau1}).
In the same vein we show that
\begin{equation}
\label{dermorceau2}
\left\|{\bf P}_{-}\left(fw_{\mu^{1/2}}(e^{\Delta_{h}}-1)\vee0\right)
\right\|^2_{L^2(\TT)}={\bf o}\left(\|h/\mu\|_{L^\infty(J)}\right).
\end{equation}
Indeed, since ${\bf P}_-$ is a contraction in $L^2(\TT)$ and $|w_{\mu^{1/2}}|\equiv1$ on $I$,
we have that
\[\left\|{\bf P}_{-}\left(fw_{\mu^{1/2}}(e^{\Delta_{h}}-1)\vee0\right)
\right\|^2_{L^2(\TT)}~\leq~
<|f|^2~,~\left|e^{\Delta_{h}}-1\right|^2>_{L^2(I)}
\]
\[
=
2<|f|^2~,~(1-\cos(\varphi/2))>_{L^2(I)}
\]
which can be treated like the right-hand side of (\ref{morceauder1}) to obtain (\ref{dermorceau2}),
granted that $|f|^2\vee0\in L\log^+L$. In view of (\ref{expcritmin}), (\ref{diffprinc}),
 (\ref{dermorceau1}) and (\ref{dermorceau2}), the proof is  complete once we have observed  that
\begin{equation}
\label{derfac}
\left\|(\mu+h)^{1/2} M\right\|_{L^2(J)}^2-
\left\|(\mu)^{1/2} M\right\|_{L^2(J)}^2=
\int_JhM^2.
\end{equation}
\hfill \boite\\

{\bf Remark:} It would be interesting to know whether Proposition \ref{cduGateau} holds true
as soon as $f\in L^2(I)$, without having to assume that $|f|^2\in L\log^+L$. 
In this case, it is easy to check using (\ref{majDini}),
(\ref{majQ2}), (\ref{morceauder1}), and  (\ref{derfac}) that
\[
\left|\Phi_M(\mu+h)-\Phi_M(\mu)-
\int_J\,h(|g_\mu|^2-M^2)\,d\theta
\right|
\]
\[
={\bf O}\left(\left(\|h/\mu\|_{L^\infty(J)}^2+\int_0^\pi\frac{\omega_{0\vee h/\mu}(t)}{t}\,dt\right)^2
\right),
\]
which is a weak substitute to (\ref{derGat}) under the (much stronger) assumption that
$0\vee h/\mu$ is Dini-continuous.

\newpage
\section{A constructive polynomial approach}
\label{const}




We now establish a finite dimensional (polynomial) analogous of Theorem  \ref{thm4} in order    
to constructively approach problem $BEP_{2, \infty}$ (\ref{eq-approx}), assuming $I$ to be a finite union of closed disjoint sub-arcs of $\TT$. This is done from the point of view of convex optimization theory  and somehow independently of the previous results, and allows us to get an alternative proof of Theorem  \ref{thm4}. 

Let $T_n$ be the space of algebraic polynomials in 
the variable $z=e^{i\theta}$ of degree less or equal to $n$ with
coefficients in $\Bbb C$. We
introduce the following finite dimensional bounded extremal problem $FBEP^n_{2, \infty}$.

\begin{description}
\item[$FBEP^n_{2, \infty}$:] For $f \in L^2(I)$, find $k_n \in T_n$ such that
$|k_n(e^{i \theta})| \leq 1$ for a.e $e^{i \theta} \in J$
and such that
\begin{equation} 
\|f - k_n\|_{L^2(I)} = \min_{\stackrel{g \in T_n}{|g| \leq 1 \,
    \mbox{a.e. in} \, J}} \|f - g\|_{L^2(I)} \, .
\label{fbepn}
\end{equation}  
\end{description} 
The existence of $k_n$ is ensured by the compactness of the
approximation set $\{g \in T_n,\, ||g||_{L^{\infty}(J)} \leq 1\}$ whereas
uniqueness follows from the convexity of this set combined with the strict
convexity of the $L^2$ norm.

For $g\in T_n$ define $$E(g)=\{x \in
\overline{J},\,\,|g(x)|=||g||_{L^{\infty}(J)}\}$$ the so called set of
critical points of $g$. 
The solution $k_n$ is  caracterised by the following result.

\begin{theorem}
\label{th_dimf}
The element $g\in T_n$ is the optimal solution of $FBEP^n_{2, \infty}$ (\ref{fbepn}) if, and only if, the following 
two conditions hold:
\begin{itemize}
\item $||g||_{L^{\infty}(J)} \leq 1$,\\
\item there exists
a set of at most $2(n+1)$ distinct points $x_i \in E(g)$ and associated positive real numbers $\lambda_i$ such that, for $r \leq 2(n+1)$:
\begin{equation}
\forall h \in T_n,\, <g-f,h>_{I}+\sum_{i=1}^{r} \lambda_i
g(x_i)\overline{ h(x_i)} =0 \, .
\label{carac_fbep}
\end{equation}
\end{itemize}
Moreover the $(\lambda_i)$ verify the following boundedness 
equation:
\begin{equation}
\label{lambda_bound}
\sum_{i=1}^{r} |\lambda_i| \leq 2||f||_{L^2(I)}^2.
\end{equation}
\end{theorem}
{\bf Remark:} 
The subset of extremal points $\{x_i,\,\,i=1,\dots ,r\}$ is possibly empty (i.e. $r=0$).\\

{\sl Proof.}
Suppose $g$ verifies the two latter conditions and differs from $k_n$. Set
$h=k_n-g$ and observe that, 
\begin{equation}
\label{neg} 
\mbox{\rm Re}\left(g(x_i)\overline{h(x_i)}\right)=\mbox{\rm Re}\left(g(x_i)\overline{ k_n(x_i)}-1\right) \leq 0
\,\,\,i=1\dots r .
\end{equation} 
From the uniqueness and optimality of $k_n$ we deduce,
\begin{equation}
\nonumber
\begin{split} 
||k_n-f||^2_{L^2(I)}&=||g-f+h||^2_{L^2(I)} \\
                  &=||g-f||^2_{L^2(I)}+||h||^2_{L^2(I)}+2\mbox{\rm Re}<g-f,h>_{I} \\
                  &<||g-f||_{L^2(I)}^2
\end{split}
\end{equation}
The latter leads to: $$ \mbox{\rm Re}<g-f,h>_{I}<0$$ which combined with
(\ref{neg}) contradicts (\ref{carac_fbep}).\\

Conversely, suppose that $g$ solves $FBEP^n_{2, \infty}$ and let $\phi_0$ be the
$\Bbb R-linear$ forms defined on $T_n$ by
$$
\forall h \in T_n, \,\,\phi_0(h)=\mbox{\rm Re} <g-f,h>_{I}.
$$
For each critical point $x$ of $g$ define the following
$\Bbb R-linear$
form $\phi_x$,
$$\forall h\in T_n,\,\, \phi_x(h)=\mbox{\rm Re}\left(g(x)\overline{h(x)}\right).$$
Finally define $K$
as follows:
$$K=\{\phi_0\} \cup \{\phi_x,\,\,x \in E(g)\}.$$ 
The set $K$ can be seen seen as a subset of the 
dual of $T_n^{\RR}$ which is defined to be the real vector space formed by
the elements of $T_n$. Simple inspection shows that $K$ is a closed
bounded set of $(T_n^{\RR})^*$, hence compact, as well as its convex hull
denoted by $\hat{K}$ (note that the finite dimensional setting is crucial 
here).
In order to get a contradiction, suppose that $0 \not \in \hat{K}$, so that
by the Hahn Banach theorem there exists $h_0 \in T_n$ such that:
$$ \forall \phi \in \hat K ,\,\, \phi(h_0) \geq \tau > 0. $$  
The latter and the continuity of $g$ and $h_0$ ensure the existence of
a neighborhood $V$ of $E(g)$ such that for $x$ in  $U=J \cap V$, then
$\mbox{\rm Re} \left(g(x)\overline{h_0(x)}\right) \geq \frac{\tau}{2}$, and for $x$ in $J
\backslash U$, then $|g(x)| \leq 1-\delta
\,\,(\delta>0) $.
Observe first that for $\epsilon >0 $ such that $\epsilon
||h_0||_{L^{\infty}(J)}< \delta$ we have 
\begin{equation}
\sup_{x \in J \backslash U} |g(x) - \epsilon h_0(x)| \leq 1. \label{surcU} 
\end{equation}
For all $x \in U$:
\begin{equation}
\nonumber
\begin{split}
|g(x)-\epsilon h_0(x)|^2 &=
|g(x)|^2-2\mbox{\rm Re} \left(g(x)\overline{h(x)}\right)+\epsilon^2|h_0(x)|^2 \\
&\leq |g(x)|^2-2\mbox{\rm Re} \left(\epsilon g(x)\overline{h(x)} \right)+\epsilon^2|h_0(x)|^2 \\
&\leq 1-\epsilon \tau+\epsilon^2||h_0||_{L^{\infty}(J)}^2 \, ,
\end{split}
\end{equation}
which combined with (\ref{surcU}) shows that for $\epsilon$ sufficiently
small we have,
\begin{equation}
||g-\epsilon h||_{L^{\infty}(J)}\leq 1 \label{surJ}.
\end{equation}
But 
\begin{equation*}
\begin{split}
||f-g-\epsilon h||_{L^2(J)}^2&=||f-g||_{L^2(J)}^2-2\epsilon
\phi_0(h)+\epsilon^2||h||_{L^2(J)}^2 \\
&\leq ||f-g||_{L^2(J)}^2-2\epsilon \tau + \epsilon^2||h||_{L^2(J)}^2 \, ,
\end{split}
\end{equation*}
which, with (\ref{surJ}), indicates that for $\epsilon$ small enough,
the function $(g-\epsilon h)$ provides a better candidate for (\ref{fbepn}) than $g$. Hence $0 \in \hat K$.

Carath\'eodory's theorem \cite{Rivlin2} is now to the effect that there exists $r'$
elements $\gamma_k$ of $K$ with $r' \leq 2(n+1)+1$ such that:
\begin{equation}
\label{convex_comb}
\sum_{i=1}^{r'} \alpha_i \gamma_i =0 \, ,\\
\end{equation}
where the $\alpha_i$ are positive real numbers such that $\sum
\alpha_i=1$.
Now suppose that $\phi_0 \neq \gamma_k$, $k = 1, \cdots, r'$. Evaluating the
sum (\ref{convex_comb}) on the element $g$ yields,
$$  \sum_{i=1}^{r'} \alpha_i \gamma_i(g)= \sum_{i=1}^{r'} \alpha_i
|g(x_i)|^2=1.$$ Equation (\ref{convex_comb}) can therefore be written
as:
$$\forall h \in T_n,\,\,\ \alpha_1 \mbox{\rm Re}<f-g,h>_{I}+\sum_{i=2}^{r'} \alpha_i
\mbox{\rm Re}(g(x_i)\overline{h(x_i)}) = 0 .$$ Dividing the
latter by $\alpha_1$ and noting that the latter equation is also true
for the element $(ih)$ leads to
(\ref{carac_fbep}).
Finally replacing $h$ by $g$ in (\ref{carac_fbep}) we obtain:
\begin{equation}
\nonumber
\begin{split}
\sum_{i=1}^{r}\lambda_i&=\sum_{i=1}^{r}|\lambda_i| =<f-g,g>_{I}\\ 
          &\leq \,\, <f-g,f-g>_{I}+|<f-g,f>_{I}|\\
          &\leq ||f||_{L^2(I)}^2+||f-g||_{L^2(I)}||f||_{L^2(I)} 
\leq 2 ||f||_{L^2(I)}^2 \, ,
\end{split}
\end{equation}
where the last two inequalities are obtained by observing that $0$ is a
valid candidate for (\ref{fbepn}). 
\hfill \boite\\
 
The next result is mainly concerned with the behavior of $FBEP^n_{2, \infty}$ when
$n$ goes towards infinity. 
\begin{theorem}
\label{conv_n}
The sequence $(k_n)$ of polynomials solutions to $FBEP^n_{2, \infty}$ (\ref{fbepn}) converges as $n \to \infty$ towards the solution $g_0 \in H^2$ to $BEP_{2, \infty}$ (\ref{eq-approx}) with
respect to the $L^2(\TT)$ norm. On $J$,  the sequence $(k_n)$  converges 
w.r.t. the weak-* topology of $L^{\infty}(J)$ and w.r.t. the $L^p(J)$ norm, for $1 \leq p < \infty$. 
In other
words we have:
$$
\lim_{n \rightarrow \infty} ||g_0 - k_n||_{L^2(\TT)}=0 \, ,
$$
$$
\forall h \in H^1,\,\,\lim_{n \rightarrow \infty}<k_n,h>_{J}=<g_0,h>_{J}  \, ,
$$
$$
\forall p \, , \ 1 \leq p < \infty, \;\;\; \lim_{n \rightarrow \infty} ||g_0 - k_n||_{L^p(J)}=0  \, . 
$$
\end{theorem}
{\sl Proof.}
Our first objective is to show (constructively, below) that $g_0$ can be approximated in the $L^2$ sense on $I$ by
polynomials that remain bounded by 1 on $J$ (thus belong to the approximating class of $FBEP^n_{2, \infty}$ for large enough $n$). 
By hypothesis, $I$ is the union of
$N$ disjoint closed sub-arcs of $\TT$ and can therefore be written as,
$$I=\bigcup_{i=1}^{N} (e^{ia_i},e^{ib_i})$$
where without loss of generality we can impose 
$$ 0=a_1 \leq b_1 \leq a_2 \dots \leq b_N \leq 2\pi.$$ 
The inner-outer decomposition of $g_0$ therefore takes the form \cite{Duren, Garnett}:
$$g_0(z)=B (z) \, \exp\left(\sum_{i=1}^{N}
\frac{1}{2\pi}\int_{a_i}^{b_i}\frac{e^{it}+z}{e^{it}-z}\log(|g_0|)dt\right).$$
Let $(\epsilon_n)$ be a decreasing sequence of positive real numbers
converging towards 0.
We define a sequence $(v_n)$ of $H^2$ as:
\begin{equation}
\nonumber
\begin{split}
v_n(z)&=B(z) \, \exp\left(\frac{1}{2\pi}\sum_{i=1}^{N}\int_{a_i+\epsilon_n}^{b_i-\epsilon_n}\frac{e^{it}+z}{e^{it}-z}\log(|g_0|)dt\right) = B (z) \times
\\
   &\exp\left(\frac{-1}{2\pi}\sum_{i=1}^{N}\int_{a_i}^{a_i+\epsilon_n}\frac{e^{it}+z}{e^{it}-z}\log(|g_0|)dt+\int_{b_i}^{b_i+\epsilon_n}\frac{e^{it}+z}{e^{it}-z}\log(|g_0|)dt\right).
\end{split}
\end{equation}
We claim that $(v_n)$ converges to $g_0$ in the $L^2$ sense on $I$.
To prove this we first show that $v_n$ converges pointwise to $g_0$
a.e. in $I$. Let $e^{i\psi}$ be a point of the
interior of $I$ and such that $g_0$ admits a radial limit at $e^{i\psi}$.
For $n$ sufficiently
large $e^{i\psi}$ is contained in none of the sub-arcs
$(a_i,a_i+\epsilon_n)$ nor $(b_i,b_i+\epsilon_n)$.
This is to the
effect that
\[
\left|\int_{a_i}^{a_i+\epsilon_n}\frac{e^{it}+e^{\psi}}{e^{it}-e^{i\psi}}\log(|g_0|)dt\right|
\leq
\int_{a_i}^{a_i+\epsilon_n}\frac{|e^{it}+e^{\psi}|}{|e^{it}-e^{i\psi}|}|\log(|g_0|)|dt\]
\[
=\int_{a_i}^{a_i+\epsilon_n}|\mbox{cotg}(\frac{t-\psi}{2})\log(|g_0|)|dt
\]
\begin{equation}
\label{maj}
\leq \max(|\mbox{cotg}(\frac{a_i-\psi}{2})| \, , \ 
|\mbox{cotg}(\frac{a_i+\epsilon_n-\psi}{2})|)\int_{a_i}^{a_i+\epsilon_n}|\log(|g_0|)|dt.
\end{equation}
The same is true with $b_i$ in place of $a_i$. As the last term of
$(\ref{maj})$ can be set arbitrarily small for $n$ sufficiently large,
the pointwise convergence of $v_n$ to $g_0$ is ensured. Finally remark that by
construction $|v_n| \leq |g_0|+1$, which leads to the majoration
$|v_n-g_0|^2 \leq (2|g_0|+1)^2$. Hence Lebesgue's dominated convergence
theorem applies and $$\lim_{n \rightarrow \infty} ||g_0-v_n||_{L^2(I)}=0.$$      

Now let $\epsilon >0$ and $0<\alpha<1$ such that $||g_0-\alpha
g_0||_{L^2(I)}\leq \frac{\epsilon}{4}$. Let $n_0$ sufficiently large such
that $||v_{n_0}-g_0||_{L^2(I)}\leq \frac{\epsilon}{4}$. For $r<1$ we
define $u_r$ belonging to disk algebra in the following way,
$$\forall \theta \in
[0,2\pi],\,\,u_r(e^{i\theta})=\int_{\TT}P_r(\theta-t)v_{n_0}(re^{it})dt,$$ 
where $P_r$ is the Poisson kernel. 

Let $e^{i\phi}
\in J$. Observe that by construction $|v_n|=1$ a.e on the sub-arc
$(e^{i(\phi-\epsilon_{n_0})},e^{i(\phi+\epsilon_{n_0}}))$. This is to the effect 
    that 
\begin{equation}
\nonumber
\begin{split}
|u_r(e^{i\phi})| & \leq  \int_{\TT}P_r(\phi-t)|v_{n_0}(re^{it})|dt \\
           & \leq  P_r(\epsilon_{n_0})\int_{\TT}|v_{n_0}(re^{it})|dt+
\int_{-\epsilon_{n_0}}^{+\epsilon_{n_0}} P_r(t)dt\\
           & \leq P_r(\epsilon_{n_0})||v_{n_0}||_{L^1(\TT)}+
          1 .
\end{split}
\end{equation}
Hence for $r$ sufficiently close to $1$, we have $|u_r| \leq 1/\alpha^2$ on
$J$ and $||u_r-v_{n_0}||^2_{L^2(I)} \leq \frac{\epsilon}{4}$. Finally,
we call $q$ the truncated Taylor expansion of $u_r$ (which converges
uniformly on $\TT$), where the truncation order has been chosen large
enough so as to ensure that $|q| \leq  1/\alpha$ on $J$ and 
$||q-u_r||^2_{L^2(I)} \leq \frac{\epsilon}{4}$. We have:
\[
||\alpha q - g_0||_{L^2(I)}
 \]
\[
\leq 
\alpha \left(||q-u_r||_{L^2(I)}+||u_r-v_{n_0}||_{L^2(I)}+
||v_{n_0}-g_0||_{L^2(I)}\right)+||g_0-\alpha g_0||_{L^2(I)}
\]
\[
\leq  \epsilon \, . 
\]
Hence, the $\alpha q$ furnish the desired polynomials.

Because they belong to the approximating class in $FBEP^n_{2, \infty}$, for large enough $n$, the above inequality is to the effect  that:
\begin{equation}
\label{eps4}
\lim_{n \rightarrow \infty}||f-k_n||_{L^2(I)}=||f-g_0||_{L^2(I)}.
\end{equation} 
As a bounded sequence of elements of $H^2$, $(k_n)$ admits a weak
convergent sub-sequence. The traces on $J$ of this sub-sequence are
bounded in the $L^{\infty}$ sense on $J$, hence up to another
sub-sequence we obtain a sequence $(k'_n)$ converging in addition in the 
weak-* sense on $J$. Let $g$ be the weak limit ($H^2$ sense) of $k'_n$. As the balls
are weak-* closed in $L^\infty$ we have $||g||_{L^{\infty}(J)}
\leq 1$, and it follows from (\ref{eps4}) that
$||f-g||_{L^2(I)}=||f-g_0||_{L^2(I)}$. The uniqueness of $g_0$ leads 
to $g=g_0$. Now (\ref{eps4}) and the constraint's saturation are to the
effect that $\limsup ||k'_n||_{L^2(\TT)} \leq ||g_0||_{L^2(\TT)}$ which 
in turn proves that $k'_n$ converges strongly in the $L^2$ sense to
$g_0$. The same kind of remark on the trace on $J$ of $k'_n$ leads to
the strong convergence in the $L^p$ (for the reflexive cases
$1<p<\infty$) sense on $J$. Finally we remark
that the preceding arguments are also true when $k_n$ is replaced by
any subsequence of the latter; hence $k_n$ contains no
sub-sequence not converging to $g_0$. \hfill \boite\\ 

{\bf Remark:} 
A discretization on $T_n$ is also taken up in \cite{schneck2}, for approximation issues of $BEP$ type in $L^p(I)$ and constrained on $\TT$, with smooth data \cite{schneck1}. This issues might themselves be normalized and formulated as $BEP_{p, \infty}$ type problems, $g$ being this time a Schur function.\\

When $I$ is a finite union of closed disjoint arcs of $\TT$, 
Theorem  \ref{thm4} may now be viewed as a corollary to Theorem \ref{th_dimf}, of which it is a infinite dimensional analogous. We  detail below this alternative proof.\\

We define $H^{2,\infty}$ and $H^{2,1}$ to be the
following vector spaces:
$$ H^{2,\infty} = \{h \in H^2,\,\,||h||_{L^{\infty}(J)}<\infty\}, $$
$$ H^{2,1}=\{h \in H^1,\,\,||h||_{L^2(I)}<\infty\}.$$
We begin with a technical lemma.
\begin{lemma}
\label{lemme_integr}
Let $v \in L^{1}(J)$ such that ${\bf P}_+(0 \vee v) \in H^{2,1}$, then the
following holds: 
$$ \forall h \in H^{2,\infty},\,\,<{\bf P}_+(
0 \vee v),h>_{\TT}=<v,h>_{J}. $$
\end{lemma}
{\sl Proof.}
Let $u$ be the function defined on $\TT$ by $$u=(0 \vee v) -{\bf P}_+(
0 \vee v).$$ By this very definition all the Fourier coefficients of $u$
of non-negative index vanish, $u$ is $L^2$ integrable on $I$ and
$L^1$ integrable on $J$. Hence we conclude that $\overline{u} \in
H^{2,1}$ and that $\overline{u}(0)=0$ ($\overline{u}$ has now a
canonical extension to the disc). Let $h \in H^{2,\infty}$.  We have:
\begin{equation}
\begin{split}
<v\chi_J,h>_{\TT}&=\,\,<u,h>_{\TT}+<{\bf P}_+(0 \vee v ),h>_{\TT} \\
                      &=\overline{u}(0)h(0)+<{\bf P}_+(0 \vee v ),h>_{\TT} \\
                      &=\,\,<{\bf P}_+(0 \vee v),h>_{\TT} \, ,
\end{split}
\end{equation}
where the second equality occurs because $(\overline{u}h) \in H^1$.
\hfill \boite\\  

{\sl Proof of Theorem  \ref{thm4}.} 
In view of (\ref{caracb-bepg}), point (ii) of Theorem \ref{thm4} and equation (\ref{carac_bepg}) can be equivalently stated as:

there exists
a non-negative function $\lambda \in L_{\RR}^1(J)$ such that, 
\begin{equation}
\forall h \in H^{2,\infty},\, <g-f,h>_{I}+<\lambda g,h>_{J}=0. 
\label{carac_bep}
\end{equation}
Suppose that $g \in H^2$ verifies $|g(e^{i\theta})|=1$ for a.e. $e^{i\theta} \in J$ and that (\ref{carac_bep}) holds, while $g \neq g_0$, the solution to $BEP_{2, \infty}$. Set $h=(g_0-g) \in H^{2,\infty} $ and observe that,
\begin{equation}
\label{neg1}
\mbox{\rm Re}<\lambda g,h>_{J}=\frac{1}{2\pi}\int_{J} \lambda(\mbox{\rm Re}(\overline{g}g_0)-1)
\leq 0 
\end{equation}
Uniqueness and optimality of $g_0$ lead (as in the proof of Theorem \ref{th_dimf})) to 
$$ \mbox{\rm Re}<g-f,h>_{I}<0 \, ,$$ which combined with
(\ref{neg1}) contradicts (\ref{carac_bep}).

Suppose now that $g$ is the optimal solution of $BEP_{2, \infty}$. The property $|g|=1$ 
on $J$ has already been proved in Theorem \ref{thm1}. In order to let $n$ go to infinity rewrite (\ref{carac_fbep}) with self explaining notations as:
\begin{equation}
\label{orth}
\forall m \in \{0 \dots n\},\,
<k_n-f,e^{im\theta}>_{I}+\sum_{i=1}^{r(n)} \lambda^n_i k_n(e^{i\theta_i^n})\overline{ e^{im\theta_i^n}} =0. 
\end{equation}
We define $(\Lambda_n)$ to be a family of linear forms
on $C(J)$ defined in the following way:
$$\forall u \in C(J),\;\;\Lambda_n(u)=\sum_{i=1}^{r(n)} \lambda^n_i
k_n(e^{i\theta_i^n})u(e^{\theta_i^n}).$$
Equation (\ref{lambda_bound}) shows now that $(\Lambda_n)$ is a bounded
sequence of elements in $C(J)^*$ which by the Banach-Alaoglu theorem \cite{Brezis}
admits a weak-* converging subsequence whose limit we call
$\Lambda$. Now the Riesz representation theorem ensures the existence of 
a complex measure $\mu$ associated to $\Lambda$, so that appealing to
Theorem \ref{conv_n} we obtain 
\begin{equation}
\forall m \in \NN,\,
<g_0-f,e^{im\theta}>_{I}+\int_{J} \overline{e^{im\theta}}d\mu =0 
\end{equation}
by taking the limit in (\ref{orth}). Now F. and M. Riesz Theorem
asserts that $\mu$ is absolutely continuous with respect to the Lebesgue measure 
so that there exists $v \in L^1(J)$ such that:
$$
\forall m \in \Bbb N,\,
<g_0-f,e^{im\theta}>_{I}+<v,e^{im\theta}>_{J} =0 \, , 
$$
which is equivalent to
\begin{equation}
\label{first_inf_equ}
\forall m \in \Bbb N,\,
<g_0-f,e^{im\theta}>_{I}+<\lambda g_0,e^{im\theta}>_{J} =0 \, ,
\end{equation}
where we have defined $\forall z \in
J,\;\;\lambda(z)=v(z)\overline{g_0(z)}$.
Equation (\ref{first_inf_equ}) is to the effect that, 
$$ {\bf P}_+((g_0-f)\chi_I)=- {\bf P}_+(0 \vee \lambda g_0)$$ 
which indicates that ${\bf P}_+(0 \vee \lambda g_0)$ is in $H^2$ (note that
this is not trivial, since $v$ occurred till now as an $L^1$
function). 
Now thanks to Lemma \ref{lemme_integr} we obtain,
\begin{equation}
\label{final_orth}
\forall u \in H^{2,\infty}  <g_0-f,u>_{I}+<\lambda g_0,u>_{J} =0.
\end{equation}
In order to prove that  $\lambda \in \RR^+$,  
consider the valid variation $h=g_0b$ where $b$ is defined as in
(\ref{def-b}),
\begin{equation}
b(z) = \frac{1}{2 \, \pi} \, \int_{I} \frac{e^{it} + 
z}{e^{it} - z} \, h(e^{it}) \, d t = \frac{1}{2 \, \pi} \,
\int_{\TT} \frac{e^{it} +  
z}{e^{it} - z} \, \chi_I(e^{it}) \, h(e^{it}) \, d t \, ,
\end{equation}
with $h \in C^\infty_{c,\RR}(I)$. We already now (as $h$ is a valid
variation) that 
$$ \mbox{\rm Re}<(f-g_0)\overline{g_0},b>_{I}=0,$$ which yields 
$$ \forall h \in  C^\infty_{c,\RR}(I),\;\;\ <\mbox{\rm Im}(\lambda),\mbox{\rm Im}(b)>_{J}=0$$ by
remarking that $b$ is pure imaginary on $J$. Now using the same
technique as in the proof of the constraint's saturation we obtain
$$ \forall u \in  C_{\RR}(J),\;\;\ <\mbox{\rm Im}(\lambda),u>_{J}=0$$ which
proves that $\lambda$ takes real values.

Finally using the fact that $BEP_{2, \infty}$ is a convex problem we obtain 
by derivating one more time that:
$$ \mbox{\rm Re}<(g_0-f),b^2>_{I}\,\, \geq 0 $$ which leads to $$ \forall u \in
C_{\RR}(J),\;\;\ <\lambda,u^2>_{J}\, \geq 0.$$ Hence $\lambda \geq 0$. Because ($\ref{final_orth}$) implies that (\ref{carac_bepg}) holds, the function $(f-g_0)\vee \lambda g_0$ cannot vanish
on a measurable set of positive measure unless it is the zero
function. But this would imply $f=g_0$ a.e on $I$ which contradicts the assumptions on $f$. This yields $\lambda>0$ a.e on $J$. \hfill \boite

\bibliographystyle{plain}
\bibliography{biblio2009}

\begin{thebibliography}{10}

\bibitem{AhernClark}
P.R. Ahern and D.N. Clark.
\newblock On functions orthogonal to invariant subspaces.
\newblock {\em Acta Math.}, 124:191--204, 1970.

\bibitem{Aizenberg}
L.~Aizenberg.
\newblock {\em Carleman's formulas in complex analysis}.
\newblock {Kluwer Academic Publishers}, 1993.

\bibitem{ablinria}
D.~Alpay, L.~Baratchart, and J.~Leblond.
\newblock Some extremal problems linked with identification from partial
  frequency data.
\newblock In J.L.~Lions R.F.~Curtain, A.~Bensoussan, editor, {\em {10th Conf.
  Analysis Optimization of Systems, Sophia-Antipolis, 1992}}, volume 185 of
  {\em LNCIS}, pages 563--573. Springer-Verlag, 1993.

\bibitem{AE98}
S.~Ansari and P.~Enflö.
\newblock Extremal vectors and invariant subspaces.
\newblock {\em Trans. Amer. Math. Soc.}, 350:539--558, 1998.

\bibitem{RRbglosw}
L.~Baratchart, J.~Grimm, J.~Leblond, M.~Olivi, F.~Seyfert, and F.~Wielonsky.
\newblock Identification d'un filtre hyperfr\'equences par approximation dans
  le domaine complexe, 1998.
\newblock INRIA technical report no. 0219.

\bibitem{BLPprep}
L.~Baratchart, J.~Grimm, J.~Leblond, and J.R. Partington.
\newblock Asymptotic estimates for interpolation and constrained approximation
  in {$H^2$} by diagonalization of {Toeplitz} operators.
\newblock {\em Integral Equations and Operator Theory}, 45:269--299, 2003.

\bibitem{partII}
L.~Baratchart and J.~Leblond.
\newblock Hardy approximation to {$L^p$} functions on subsets of the circle
  with $1 \leq p < \infty$.
\newblock {\em Constructive Approximation}, 14:41--56, 1998.

\bibitem{partI}
L.~Baratchart, J.~Leblond, and J.R. Partington.
\newblock Hardy approximation to {$L^\infty$} functions on subsets of the
  circle.
\newblock {\em Constructive Approximation}, 12:423--436, 1996.

\bibitem{aakblp}
L.~Baratchart, J.~Leblond, and J.R. Partington.
\newblock Problems of {Adamjan--Arov--Krein} type on subsets of the circle and
  minimal norm extensions.
\newblock {\em Constructive Approximation}, 16:333--357, 2000.

\bibitem{BS02}
L.~Baratchart and F.~Seyfert.
\newblock An {$L^p$} analog to {AAK} theory for $p \geq 2$.
\newblock {\em J. Funct. Anal.}, 191, 2002.

\bibitem{BeSh}
C.~Bennett and R.~Sharpley.
\newblock {\em Interpolation of operators}.
\newblock Number 129 in Pure and Applied Mathematics. Academic Press, 1988.

\bibitem{BorLew}
J.M. Borwein and A.S. Lewis.
\newblock {\em Convex Analysis and Nonlinear Optimization}.
\newblock {CMS} Books in {M}ath. Can. Math. Soc., 2006.

\bibitem{Brezis}
H.~Br\'ezis.
\newblock {\em Analyse fonctionnelle}.
\newblock {Dunod}, 1999.

\bibitem{CP2002}
I.~Chalendar and J.~R. Partington.
\newblock Constrained approximation and invariant subspaces.
\newblock {\em J. Math. Anal. Appl.}, 289(1):176--187, 2003.

\bibitem{CPM2002}
I.~Chalendar, J.~R. Partington, and M.~P. Smith.
\newblock Approximation in reflexive {Banach} spaces and applications to the
  invariant subspace problem.
\newblock {\em Proc. A.M.S.}, 132(4):1133--1142, 2004.

\bibitem{DFT}
J.~C. Doyle, B.~A. Francis, and A.~R. Tannenbaum.
\newblock {\em Feedback Control Theory}.
\newblock Macmillan Publishing Company, 1992.

\bibitem{Duren}
P.L. Duren.
\newblock {\em Theory of $H^p$ spaces}.
\newblock {Academic Press}, 1970.

\bibitem{Fuhrmann}
P.A. Fuhrmann.
\newblock {\em Linear systems and operators in Hilbert spaces}.
\newblock {McGraw--Hill}, 1981.

\bibitem{Garnett}
J.B. Garnett.
\newblock {\em Bounded analytic functions}.
\newblock {Academic Press}, 1981.

\bibitem{GoKry}
G.~M. Goluzin and V.~I. Krylov.
\newblock Generalized carleman formula and its application to analytic
  extension of functions.
\newblock {\em Mat. Sb}, 40, 1933.

\bibitem{Isakov}
V.~Isakov.
\newblock {\em {Inverse problems for partial differential equations}}.
\newblock Number 127 in {Applied Mathematic Sciences}. {Springer}, 1998.

\bibitem{imh2}
B.~Jacob, J.~Leblond, J.-P. Marmorat, and J.R. Partington.
\newblock A constrained approximation problem arising in parameter
  identification.
\newblock {\em Linear Algebra and its Applications}, 351-352:487--500, 2002.

\bibitem{Koosis}
P.~Koosis.
\newblock {\em Introduction to {$H_p$} spaces}.
\newblock {Cambridge University Press}, 1980.

\bibitem{KN}
M.G. Krein and P.Y. Nudel'man.
\newblock {Approximation of $L^2(\omega_1, \omega_2)$ functions by minimum--
  energy transfer functions of linear systems}.
\newblock {\em {Problemy Peredachi Informatsii}}, 11(2):37--60, 1975.
\newblock English translation.

\bibitem{Yosida}
K.Yosida.
\newblock {\em Functional Analysis}.
\newblock Grundlehren der Math. Wissenschaften. {Springer--Verlag}, 1980.

\bibitem{Lavrentiev}
M.M. Lavrentiev.
\newblock {\em Some improperly posed problems of mathematical physics},
  volume~11 of {\em tracts in {Natural Philosophy}}.
\newblock {Springer}, 1967.

\bibitem{LMP08}
J.~Leblond, J.-P. Marmorat, and J.R. Partington.
\newblock Solution of inverse diffusion problems by analytic approximation with
  real constraints.
\newblock {\em J. Inv. Ill-Posed Problems}, 16(1):89--105, 2008.

\bibitem{LO98}
J.~Leblond and M.~Olivi.
\newblock Weighted {$H^2$} approximation of transfer functions.
\newblock {\em MCSS (Math. Control, Signals, Systems)}, 11:28--39, 1998.

\bibitem{LP99}
J.~Leblond and J.~R. Partington.
\newblock Constrained approximation and interpolation in hilbert function
  spaces.
\newblock {\em J. Math. Anal. Appl.}, 234(2):500--513, 1999.

\bibitem{Nikolskii2}
N.~K. Nikolskii.
\newblock {\em Operators, functions, and systems: an easy reading}, volume~92
  of {\em Mathematical surveys and monographs}.
\newblock Amer. Math. Soc., 2002.

\bibitem{Par97}
J.R. Partington.
\newblock {\em {Interpolation, identification and sampling}}.
\newblock {Oxford University Press}, 1997.

\bibitem{Patil}
D.J. Patil.
\newblock Representation of {$H^p$} functions.
\newblock {\em {Bull. A.M.S.}}, 78(4):617--620, 1972.

\bibitem{Pichorides}
S.~K. Pichorides.
\newblock On the best values of the constants in the theorems of m. riesz,
  zygmund and kolmogorov.
\newblock {\em {Studia Math.}}, 44:165--179, 1972.

\bibitem{Rivlin2}
T.J. Rivlin.
\newblock {\em Chebyshev polynomials}.
\newblock {Wiley-Interscience}, 1990.

\bibitem{Rudin}
W.~Rudin.
\newblock {\em Real and complex analysis}.
\newblock {McGraw--Hill}, 1987.

\bibitem{schneck1}
A.~Schneck.
\newblock Constrained optimization in hardy spaces.
\newblock Preprint, 2009.

\bibitem{schneck2}
A.~Schneck.
\newblock Constrained optimization in hardy spaces {II}: Numerics.
\newblock Preprint, 2009.

\bibitem{thesefab}
F.~Seyfert.
\newblock {Probl\`emes extr\'emaux dans les espaces de Hardy}.
\newblock These de Doctorat, Ecole des Mines de Paris, 1998.

\bibitem{MSmith}
M.~Smith.
\newblock Constrained approximation in banach spaces.
\newblock {\em Constructive Approximation}, pages 465--476, 2003.

\end{thebibliography}

\newpage
\tableofcontents

\end{document}